\documentclass[12pt]{article}

\usepackage{graphicx}
\usepackage{psfrag}
\usepackage{amsmath,amssymb,mathrsfs}
\usepackage{mathptmx}
\usepackage{ntheorem}
\usepackage{graphicx,mdframed}

\makeatletter
\newsavebox\myboxA
\newsavebox\myboxB
\newlength\mylenA
\newcommand*\xoverline[2][0.75]{%
    \sbox{\myboxA}{$\m@th#2$}%
    \setbox\myboxB\null
    \ht\myboxB=\ht\myboxA%
    \dp\myboxB=\dp\myboxA%
    \wd\myboxB=#1\wd\myboxA
    \sbox\myboxB{$\m@th\overline{\copy\myboxB}$}
    \setlength\mylenA{\the\wd\myboxA}
    \addtolength\mylenA{-\the\wd\myboxB}%
    \ifdim\wd\myboxB<\wd\myboxA%
       \rlap{\hskip 0.5\mylenA\usebox\myboxB}{\usebox\myboxA}%
    \else
        \hskip -0.5\mylenA\rlap{\usebox\myboxA}{\hskip 0.5\mylenA\usebox\myboxB}%
    \fi}
\makeatother

\font\tencmmib=cmmib10 \skewchar\tencmmib '60
\newfam\cmmibfam
\textfont\cmmibfam=\tencmmib

\def\lessim{\ \lower4pt\hbox{$
		\buildrel{\displaystyle <}\over\sim$}\ }
\def\gessim{\ \lower4pt\hbox{$\buildrel{\displaystyle >}
		\over\sim$}\ }

\def\si{{\sigma}}

\def\qe{{q_{\mbox{\tiny\rm EA}}}}

\def\eps{{\varepsilon}}
\def\ch{{\mbox{\rm ch}\hspace{0.4mm}}}

\def\tap{\mbox{\tiny{\rm TAP}}}
\def\qp{q_{\mbox{\tiny{\rm P}}}}
\def\ap{\alpha_{\mbox{\tiny{\rm P}}}}

\def\qed{\hfill\hbox{\rlap{$\sqcap$}$\sqcup$}}

\newcommand{\la}{\langle}
\newcommand{\ra}{\rangle}

\newcommand{\e}{\mathbb{E}}
\newcommand{\p}{\mathbb{P}}
\newcommand{\Reals}{\mathbb{R}}

\newcommand{\bsigma}{{\xoverline{\sigma}}}
\newcommand{\brho}{{\xoverline{\rho}}}

\newcommand{\PP}{{\cal P}}

\newcommand{\myth}{{\mbox{\rm th}}}

\newtheorem{lemma}{\bf Lemma}

\newtheorem{theorem}{\bf Theorem}

\newtheorem{remark}{\bf Remark}

\newtheorem{proposition}{\bf Proposition}

\newenvironment{Proof of lemma}{\noindent{\bf Proof of Lemma}}{\hfill$\Box$\newline}
\newenvironment{Proof of theorem}{\noindent{\it Proof of Theorem}}{\hfill\scriptsize{$\Box$}\newline}
\newenvironment{Proof of theorems}{\noindent{\bf Proof of Theorems}}{\hfill$\Box$\newline}
\newenvironment{Proof of proposition}{\noindent{\bf Proof of Proposition}}{\hfill$\Box$\newline}
\newenvironment{Proof of propositions}{\noindent{\bf Proof of Propositions}}{\hfill$\Box$\newline}
\newenvironment{Proof of exercise}{\noindent{\it Proof of Exercise:}}{\hfill$\Box$}
\begin{document}
	
		\title{On the TAP free energy in the mixed $p$-spin models}
		
		\author{
			Wei-Kuo Chen\thanks{School of Mathematics, University of Minnesota. Email: wkchen@umn.edu}
		    \and
			Dmitry Panchenko\thanks{Department of Mathematics. University of Toronto. Email: panchenk@math.toronto.edu}
		}
		\maketitle
		
		\begin{abstract}		
			In \cite{TAP}, Thouless, Anderson, and Palmer derived a representation for the free energy of the Sherrington-Kirkpatrick model, called the TAP free energy, written as the difference of the energy and entropy on the extended configuration space of local magnetizations with an Onsager correction term. In the setting of mixed $p$-spin models with Ising spins, we prove that the free energy can indeed be written as the supremum of the TAP free energy over the space of local magnetizations whose Edwards-Anderson order parameter (self-overlap) is  to the right of the support of the Parisi measure. Furthermore, for generic mixed $p$-spin models, we prove that the free energy is equal to the TAP free energy evaluated on the local magnetization of any pure state.
		\end{abstract}
		
		\section{Introduction and main results}\label{sec1}
		
		The study of the Sherrington-Kirkpatrick (SK) model \cite{SK} using replica method culminated in the celebrated Parisi ansatz \cite{Parisi79, Parisi, Parisi83}, which over time was developed and understood much better in the physics literature (we refer to the classical account in M\'{e}zard-Parisi-Virasoro \cite{MPV}). Since then, the model has received a lot of attention in the probability and mathematical physics communities, with many rigorous results obtained in recent years. 
		
		Different from the replica method, the approach of Thouless, Anderson and Palmer in \cite{TAP} computed the free energy of the SK model by diagrammatic expansion of the partition function, deriving a representation (known as the TAP free energy) in terms of the energy, classical entropy for product measure with the same marginals as the Gibbs measure, and Onsager correction term, all defined on the expanded configuration space of local magnetizations. Their approach involved convergence conditions, which are certain constraints on the so called Edwards-Anderson (EA) order parameter of the model similar to the one that will appear in our first main result, Theorem \ref{thm2} below. The fact that the TAP equations can be used only in a region of Edwards-Anderson parameter has been also argued, for example, in Anderson \cite{PWA}. The general convergence criterion intended for the whole temperature range was later derived by Plefka in \cite{Plefka}. In addition to the TAP free energy, \cite{TAP} also described the famous TAP equations for local magnetizations, which are just the critical point equations for the TAP free energy. The connection between the TAP approach and Parisi ansatz has been studied in the physics literature (see e.g. \cite{Betal, CGPM,DDYo}), with the general idea  that (some) critical points of the TAP free energy correspond to the pure states in the Parisi ansatz. Our second main result, Theorem \ref{ThHTtapLT}, shows (when combined with Theorem \ref{thm2}) that pure states are indeed maximizers of the TAP free energy in the aforementioned region of the Edwards-Anderson order parameter, although this does not immediately translate into a statement about critical points, because the result holds in the thermodynamic limit. 
	
		Beyond the SK model, the framework of the TAP approach has  been used widely in various systems with huge complexity. These include, for instance, the random optimization problems arising from neural computation and computer science. A thoroughly discussion along these lines can be found in the workshop proceedings edited by Opper-Sadd \cite{OS}. 
		
		On the mathematics side, the TAP equations for local magnetizations were first proved by Talagrand \cite{Talbook1}, and also Chatterjee \cite{chatt10}, in the SK model at high temperature. More recently, an iterative scheme for constructing the solutions for the system of these TAP equations was introduced by Bolthausen in \cite{EB}, and was shown to converge in the entire high temperature region. Even more recently, the TAP equations for local magnetization were derived by Auffinger-Jagannath \cite{AJ-tap} for the generic mixed $p$-spin models at low temperature (in a slightly weaker sense than earlier results at high temperature, which showed that the equations hold for all spins simultaneously with high probability, while \cite{AJ-tap} proves that they hold for most of the spins with high probability). The aim of this paper is to investigate the TAP free energy representation rather than its critical points, in the setting of mixed $p$-spin models with Ising spins. 
		
		Our first main result in Theorem \ref{thm2} below states that the free energy can be written as the supremum of the TAP free energy if the Edwards-Anderson parameter is constrained to the right of the support of the Parisi measure (the functional order parameter encoding the distribution of one overlap). In a different direction we show that, if one samples a random spin configuration from the Gibbs measure and considers a pure state around it, with high probability the TAP free energy evaluated at the barycentre (local magnetization) of this pure state approximates the free energy of the whole system. In Theorem \ref{ThHTtap} we prove this for arbitrary mixed $p$-spin model at high temperature, in which case the entire system is in a pure state, and in Theorem \ref{ThHTtapLT} we prove this for generic models at low temperature. An interesting byproduct of our calculation stated in Theorem \ref{ThHTent} shows that, at high temperature, the entropy of the Gibbs measure is different from the entropy of the product measure with the same marginals, despite a well-known fact that the two measures are close on finite sets of coordinates. This gives a negative answer to the Conjecture 1.4.18 in \cite{Talbook1}. We remark that the same discrepancy holds at low temperature for the entropy of a pure state.
		
	We should mention that all our results are proved by using, essentially, the entire arsenal of mathematical theory of the Parisi solution of these models, by calculating the TAP free energy under the constraint mentioned above and comparing it with the Parisi formula for the free energy. It would be of great interest to find a direct rigorous proof relating these two formulas, in which case one would need to use a different constraint that does not refer to any properties of the Parisi solution, perhaps, in terms of the Plefka condition (\ref{Plefka}) below.
				
		\subsection{TAP free energy representation}
		
	For any integer $N\geq 1,$ consider the hypercube $\Sigma_N=\{-1,+1\}^N.$ The Hamiltonian of the mixed $p$-spin model is given by
		\begin{align}\label{eq3}
		H_N(\sigma)&=X_N(\sigma)+h\sum_{i=1}^N\sigma_i
		\end{align}
		for $\sigma=(\sigma,\ldots,\sigma_N)\in \Sigma_N,$ 
		where $X_N$ is a Gaussian process defined as
		\begin{align}
		\label{hamx}
		X_N(\sigma)=\sum_{p\geq 2}\frac{\beta_p}{N^{(p-1)/2}}\sum_{1\leq i_1,\ldots,i_p\leq N}g_{i_1,\ldots,i_p}\sigma_{i_1}\cdots\sigma_{i_p}.
		\end{align}
		Here, $g_{i_1,\ldots,i_p}$ are i.i.d. standard Gaussian random variables and $(\beta_p)_{p\geq 2}$ is a real sequence with $\sum_{p\geq 2}2^p\beta_p^2<\infty.$ Under this assumption, the series in (\ref{hamx}) is well defined since the covariance
		\begin{align}
		\e X_N(\sigma^1)X_N(\sigma^1)=N\xi\bigl(R(\sigma^1,\sigma^2)\bigr),
		\end{align} 
		where $
		\xi(s):=\sum_{p\geq 2}\beta_p^2s^p
		$
		and
		$R(\sigma^1,\sigma^2):=N^{-1}\sum_{i=1}^N\sigma_i^1\sigma_i^2$ is the overlap between $\sigma^1$ and $\sigma^2.$ An important quantity associated to the mixed $p$-spin model is the free energy,
		\begin{align}
		F_N&=\frac{1}{N}\log \sum_{\sigma\in \Sigma_N}\exp H_N(\si).
		\end{align}
		It is now well known that the limit of the free energy exists almost surely and is given by the Parisi formula. More precisely, let $\mathcal{M}$ be the space of all cumulative distribution functions on $[0,1]$ equipped with the $L_1(ds)$ distance. For any $\alpha\in \mathcal{M},$ let $\Phi_{ \alpha}(0,h)$ be the unique weak solution (see \cite{JT16}) of nonlinear partial differential equation
		\begin{align*}
		\partial_s\Phi_{ \alpha}=-\frac{ \xi''}{2}\bigl(\partial_{xx}\Phi_{ \alpha}+\alpha(s)(\partial_x\Phi_{ \alpha})^2\bigr),\,\,(s,x)\in[0,1)\times \mathbb{R}
		\end{align*}
		with the boundary condition $\Phi_{ \alpha}(1,x)=\log\mathrm{ch} (x).$
		Define a functional $\mathcal{P}$ on $\mathcal{M}$ by
		\begin{align}\label{pf}
		\mathcal{P}(\alpha)&=\log 2+\Phi_{ \alpha}(0,h)-\frac{ 1}{2}\int_0^1\xi''(s)s\alpha(s)\, ds.
		\end{align}
		The Parisi formula for the limiting free energy states that
		\begin{align}\label{parisi}
		F:=\lim_{N\rightarrow\infty}\e F_N=\inf_{\alpha\in \mathcal{M}}\mathcal{P}(\alpha).
		\end{align}
		The validity of this formula was first proved by Talagrand \cite{Tal03} for mixed even $p$-spin models and was later extended to general mixtures in Panchenko \cite{Pan00}. It was shown by Guerra \cite{Guerra} that $\mathcal{P}$ defines a Lipschitz functional on $\mathcal{M}.$ In Auffinger-Chen \cite{AC14}, it was further investigated that the functional $\mathcal{P}$ is strictly convex on $\mathcal{M},$ so the optimization problem in \eqref{parisi} has unique minimizer. Throughout this paper, this minimizer will be denoted by ${\ap}$ and we call ${\ap}(ds)$ the Parisi measure. We denote by $\ap(A)$ the measure of $A$ with respect to ${\ap}(ds)$ if $A$ is a measurable subset of $[0,1].$ We denote the largest point in the support of ${\ap}(ds)$ by
\begin{equation*}
\qp : = \max \mathrm{supp}(\ap).
\end{equation*}

To motivate the TAP free energy representation, let us recall that the classical Gibbs variational principle advocates that 
		\begin{align}\label{FGibbsRep}
		F_N&=\max_{\mu}\Bigl(\frac{\e_\mu H_N(\sigma)}{N}-\frac{E(\mu)}{N}\Bigr),
		\end{align}
		where the maximum is taken over all probability measures $\mu$ on $\Sigma_N.$ Here $\e_\mu$ is the expectation with respect to the probability measure $\mu$ on $\Sigma_N$ and $E(\mu)$ is the entropy of $\mu$, i.e., $$
		E(\mu):=\sum_{\sigma\in\Sigma_N}\mu(\sigma)\log \mu(\sigma).
		$$
		The maximum is attained by the Gibbs measure 
		\begin{equation}
		\label{GibbsGN}
		G_N(\sigma)=\frac{\exp {H_N(\sigma)}}{\sum_{\sigma'}\exp H_N(\sigma')}.
		\end{equation}
		(As usual, we will denote by $\la\,\cdot\,\ra$ the average with respect to the Gibbs measure.)
		In particular, if $\mu$ is a product measure with $m_i=\e_\mu\sigma_i$  then 
		$
		\mu(\sigma_i=\pm 1)={(1\pm m_i)}/{2}
		$
		and
		\begin{align*}
		\e_\mu H_N(\sigma)-E(\mu)&=H_N(m)-\sum_{i=1}^NI(m_i),
		\end{align*}
		where $I(x)$ is the entropy of a Bernoulli random variable on $\{-1,1\}$ with mean $x$, 
		\begin{align}\label{eq-7}
		I(x)&:=\frac{1+x}{2}\log\frac{1+x}{2}+\frac{1-x}{2}\log \frac{1-x}{2},\,\,x\in[-1,1].
		\end{align}
		A natural question is whether $F_N$ can be written as the supremum of $\e_\mu H_N(\sigma)-E(\mu)$ over product measures in the thermodynamic limit $N\rightarrow\infty$, which happens in various examples covered by the theory of non-linear large deviations in \cite{ChDem} (see also \cite{Eldan1, Eldan2}). For spin glass models, this turns out to be not the case and, in fact, the derivation of Thouless, Anderson, and Palmer \cite{TAP} produced a correction term. As we will see in the next section, the correction term is needed even at high temperature when any finite set of coordinates is asymptotically uncorrelated.
		
		The TAP correction term is expressed in terms of the function
		\begin{align}\label{eq-8}
		C(u)&:=\frac{1}{2}\bigl(\xi(1)-\xi(u)-\xi'(u)(1-u)\bigr),\,\,u\in[0,1].
		\end{align}
		The convexity of $\xi$ on $[0,1]$ implies that $C(u)\geq 0$ for $u\in [0,1]$. (The original paper \cite{TAP} only dealt with the SK model corresponding to $\xi(s)=\beta^2s^2/2$, in which case $C(u)=\beta^2(1-u)^2/4.$) The TAP free energy is then given by
		\begin{align}\label{tap}
		F_{N}^{\tap}(m)&=\frac{H_N(m)}{N}-I_N(m)+C(\qe),\,\,m\in [-1,1]^N,
		\end{align}
		where $I_N(m)=N^{-1}\sum_{i=1}^NI(m_i)$, and where the self-overlap
		\begin{align}
		\label{qeaA}
		\qe&:=R(m,m)=\frac{1}{N}\sum_{i=1}^Nm_i^2
		\end{align}
		is called the  Edwards-Anderson order parameter. The free energy is then approximated by the maximum of the TAP free energy over configurations satisfying some constraints to ensure the convergence of the diagrammatic expansion. Our result here shows that this is indeed the case under the constraint $\qe\geq \qp.$
		\begin{theorem}
			\label{thm2}
			For any mixture parameter $\xi$ and external field $h,$ 
			\begin{align}\label{thm2:eq1}
			F&=\lim_{\varepsilon\downarrow 0}\lim_{N\rightarrow\infty}\e\max_{m\in [-1,1]^N:\qe \in [{\qp}-\varepsilon,1]}F_N^{\tap}(m).
			\end{align}
			Moreover, the maximum is attained at $\qp$, so (\ref{thm2:eq1}) also holds with the constraint $\qe \in [{\qp}-\eps,\qp+\eps].$
		\end{theorem}
		
		In the classical SK model, when $\xi(s)=\beta^2s^2/2$, it was predicted by Plefka \cite{Plefka} that $F$ can be expressed as a Boltzmann entropic principle using the TAP free energy,
		$
		F=\lim_{N\rightarrow\infty}\max_{m} F_{N}^{\tap}(m),
		$
		where the maximum is taken over all $m\in[-1,1]^N$ that satisfy the  condition,
		\begin{equation}\label{Plefka}
		1>\beta^2\bigl(1-2\qe+r\bigr)
		\end{equation}
		for $r:=N^{-1}\sum_{i=1}^Nm_i^4.$ This condition does not explicitly refer to the Parisi formula and, as we mentioned above, it would be interesting to prove this representation directly.
		
		The main ingredient in the proof of Theorem \ref{thm2} is a subtle analysis of the maximum energy $$\lim_{\varepsilon\downarrow 0}\lim_{N\rightarrow\infty}\max_{m\in[-1,1]^N:|\qe-u|\leq \varepsilon}\Bigl(\frac{H_N(m)}{N}-I_N(m)\Bigr).$$
		We first establish in full generality the analogue of the Parisi variational formula for this limit on the space of the functional order parameters given by the cumulative distribution functions $\gamma(s)$ for $s\in [0,u]$ induced by positive measures on $[0,u]$ (see Theorem \ref{thm1} below). From this, by the virtue of the entropy term $I_N,$ our Parisi functional is essentially different from the original Parisi functional $\mathcal{P}$ only by the correction term $C(u)$. This gives a lower bound for the Parisi formula \eqref{parisi} as an immediate consequence of the fact that the Parisi formulas are minimization problems. More importantly, when $u=\qp,$ we show that $\gamma(s):={\ap}(s)$ for $s\in[0,\qp]$ is the minimizer to our Parisi formula of the maximum energy and consequently, the free energy equals the TAP free energy when $u=\qp$ in the thermodynamic limit.  
\subsection{TAP free energy for pure states}

One expects the TAP free energy to be approximately equal to the free energy, $F_N \approx F_N^{\tap}(m),$ when the vector $m$ is `a magnetization of a pure state'. 

\subsubsection{High temperature case}
We will discuss more what `pure state' means at low temperature, but first let us consider the simpler case of high temperature when the overlap concentrates near some constant $q\in[0,1]$,
\begin{equation}
\lim_{N\to\infty}\e\la(R_{1,2}-q)^2\ra = 0.
\label{HTcond}
\end{equation}
For example, for pure $2$-spin model corresponding to $\xi(s)=\beta^2s^2/2,$ $q$ is defined as the unique solution of $q = \e \myth^2(z\sqrt{\xi'(q)}+h)$ for $z$ a standard Gaussian random variable when $h\not =0$ and $q=0$ when $h=0,$ and (\ref{HTcond}) is known to hold for $\beta<1/2$ (see Theorem 1.4.1 in \cite{Talbook1}). The condition \eqref{HTcond} ensures (see Proposition 1.6.8 in \cite{Talbook1} and Theorem 1 in \cite{AC14}) that the Parisi measure is a Dirac measure $\ap(dt)=1_{\{q\}}(t).$  In this case, the entire system is said to be in a pure state and we take $m=(\la \sigma_1\ra,\ldots, \la\sigma_N\ra)$ to be the barycentre of the Gibbs measure $G_N.$ Recall the TAP free energy $F_N^{\tap}$ from \eqref{tap}. Note that $\qe$ associated to $m$ is approximately $q$ by \eqref{HTcond}. The following holds.
\begin{theorem} \label{ThHTtap}
	If (\ref{HTcond}) holds for some $q\in[0,1]$ and $m=(\la \sigma_1\ra,\ldots, \la\sigma_N\ra)$ then
	\begin{equation*}
	\lim_{N\to\infty}\e\bigl(F_N^{\tap}(m) -F_N\bigr)^2 = 0.
	\end{equation*}
\end{theorem}
We will prove this by showing that each term in $F_N^{\tap}(m)$ concentrates (it is well known that the free energy $F_N$ concentrates) and computing their expected values. 

Although the Gibbs representation of the free energy (\ref{FGibbsRep}) will not be used in any way in the proof of Theorem \ref{ThHTtap}, it is interesting how various terms in $F_N^{\tap}(m)$ are related to this representation. The term $-I_N(m)$ in $F_N^{\tap}$ is the entropy of the product measure on $\{-1,+1\}^N$ with the same marginals as $G_N$, with the means $m_i=\la \sigma_i\ra.$ A well known result, Proposition 1.4.14 in \cite{Talbook1}, states that at high temperature any finite number of coordinates $\sigma_i$ are asymptotically independent, so one might expect that the entropy $-\la \log G_N(\sigma)\ra/N$ is approximated by $-I_N(m).$ However, as we will see below, this turns out not to be the case and one needs to add the TAP correction term to relate these two entropies. For example, for pure SK model the correction between the two entropies is given as follows.
\begin{theorem}\label{ThHTent}
	If $\xi(s)=\beta^2 s^2/2$ and (\ref{HTcond}) holds, then 
	\begin{equation}\label{EntCorr}
	\frac{1}{N}\la \log G_N(\sigma)\ra \approx I_N(\la\sigma\ra) + C(q).
	\end{equation}
\end{theorem}
We remark that the same formula for the discrepancy of the entropy of a pure state (defined in the next section) also holds at low temperature, which can be seen from Lemma \ref{lem6} below. 

Theorem \ref{ThHTent} combined with Lemma 4.4 in \cite{Austin} implies that the Wasserstein distance with the Hamming cost function $d(\sigma^1,\sigma^2)=N^{-1}\sum_{i\leq N}\mathrm{I}(\sigma_i^1\not = \sigma_i^2)$ between $G_N$ and the product measure with the same marginals does not go to zero even when (\ref{HTcond}) holds. This gives a negative answer to the Conjecture 1.4.18 in \cite{Talbook1}.

Note that the Gibbs measure optimizes the Gibbs variational principle \eqref{FGibbsRep}. It is important to point out that the term $C(q)=C(\qp)$ in (\ref{EntCorr}) is not the only correction responsible for the term $C(\qp)$ in the TAP representation (\ref{tap}). In fact, notice that by itself it would result in the correction in the \emph{opposite} direction! This means that another correction $2C(\qe)$ is coming from the discrepancy between $\la H_N(\sigma)/N\ra$ and $ H_N(\la\sigma\ra)/N.$

\subsubsection{Low temperature case}

To study the TAP free energy within `pure states' at low temperature rigorously, we will work with generic mixed $p$-spin models defined in Section 3.7 in \cite{Pan}.

\medskip
\noindent \textbf{Definition.} \emph{(Generic $p$-spin model)} We will call the mixed $p$-spin Hamiltonian (\ref{hamx}) \emph{generic} if linear span of constants and functions $x^p$ for $p\geq 2$ such that $\beta_p\not = 0$ is dense in $C([-1,1],\|\,\cdot\,\|_\infty)$.

\medskip
This is a particularly nice subclass of $p$-spin models, because, in this case, the distribution of the overlap $R_{1,2}$ is known to converge to the Parisi measure $\ap$ and the (unique) limiting distribution of the overlap array $(R_{\ell,\ell'})_{1\leq \ell,\ell'}$ is determined by $\ap$ via the Parisi ultrametric ansatz (see Section~3.7 in \cite{Pan}). Moreover, as will be important to us, the distribution of spins over replicas $(\sigma_i^{\ell})_{i,\ell\geq 1}$ is asymptotically determined by $\ap$ in a way explained in Chapter 4 in \cite{Pan}. In fact, the distribution of spins was determined in \cite{Pan} under regularizing perturbations, but it was observed in \cite{AJ-tap} that for generic models the proof works without perturbations.

What are pure states at low temperature when the condition (\ref{HTcond}) fails? First of all, the notion of pure states is well understood in the thermodynamic limit via the so-called asymptotic Gibbs measure $G$. If $\qp $ in the largest point in the support of the Parisi measure $\ap$ then $G$ is a random measure on the sphere of radius $\sqrt{\qp }$ on a Hilbert space such that the scalar products of the i.i.d. sample from this measure have the distribution of the overlap array $(R_{\ell,\ell'})_{1\leq \ell,\ell'}$ in the limit (see, e.g., \cite{Pan}; this notion was originally defined in \cite{AA2007} using the Dovbysh-Sudakov representation). Points in the support of this measure $G$ correspond to the physicists' pure states. Moreover, under the condition $\ap(\{\qp \}) > 0$, the measure $G$ is supported by countably many points (pure states) each carrying positive random weight (see Lemma 2.7 in \cite{Pan}). 

It was explained in Jagannath \cite{JAU} how the ultrametric structure of the asymptotic Gibbs measure $G$ can be used to define approximate pure states for finite $N,$ which are clusters of spin configurations on $\{-1,+1\}^N$ that satisfy various natural properties, most importantly, the analogue of (\ref{HTcond}). Namely, with respect to the conditional Gibbs measure on a given cluster, the overlap $R_{1,2}$ concentrates near $\qp .$ How can one discover these clusters, or pure states of $G_N$? The construction turns out to be very natural, although non-deterministic. Just like points in the support of $G$ can be discovered by sampling from $G$, if we sample a random point $\sigma$ from $G_N$ then its neighbourhood 
\begin{equation}\label{PureState}
\Sigma_N(\sigma) =\bigl\{\rho \in \{-1,+1\}^N \,:\, R(\sigma,\rho)\geq \qp -\eps\bigr\}
\end{equation}
for $\eps=\eps_N \downarrow 0$ slowly enough can be viewed as one of these approximate pure states. Moreover, if we consider a large sample $\sigma^1,\ldots,\sigma^n \sim G_N$ then the clusters $\Sigma_N(\sigma^\ell)$ either almost coincide or are almost disjoint in the sense of measure $G_N$, so these random clusters can be used to decompose the Gibbs measure into approximate pure states (see \cite{JAU} for details). 

As a result, we will take the following point of view. We will treat the random set (\ref{PureState}) as a realization of a pure state, and we will show that the TAP free energy representation holds with large probability over the choice of $\sigma$, as $\eps\downarrow 0$, which means that it holds for a typical pure state. Let us state the result precisely. Consider the Gibbs weight of this random pure state, 
$$
W_N(\sigma) = G_N(\Sigma_N(\sigma)) = G_N(\rho \,:\, R(\sigma,\rho)\geq \qp -\eps),
$$
where we keep the dependence of $W_N(\sigma)$ on $\eps$ implicit. For $\sigma\in \{-1,+1\}^N,$ let $\la\,\cdot\,\ra_\sigma$ denote the average with respect to the conditional Gibbs measure $G_N$ on $\Sigma_N(\sigma),$
\begin{equation*}
\la f(\rho^1,\ldots,\rho^n) \ra_\sigma := \frac{1}{W_N(\sigma)^n} \sum_{\rho^1,\ldots,\rho^n} f(\rho^1,\ldots,\rho^n)G_N(\rho^1)\cdots G_N(\rho^n).
\end{equation*}
We will denote by $m(\sigma)$ the barycentre of $\Sigma_N(\sigma)$ (magnetization of $G_N$ on $\Sigma_N(\sigma)$),
\begin{equation}\label{emdef}
m(\sigma)=(m_1(\sigma),\ldots, m_N(\sigma))=( \la \rho_1 \ra_\sigma,\ldots,  \la \rho_N \ra_\sigma).
\end{equation}
The following holds.
\begin{theorem} \label{ThHTtapLT}
	For generic mixed $p$-spin models, if $m(\sigma)$ is defined in (\ref{emdef}) then
	\begin{equation}\label{ThHTtapLT:eq1}
	\lim_{\eps\downarrow 0}\limsup_{N\to\infty}
	\e\bigl\la \bigl( F_N^{\tap}(m(\sigma)) - F_N \bigr)^2\bigr\ra = 0.
	\end{equation}
\end{theorem}
(Note that here we do not make any assumptions on $\ap$, for example, $\ap(\{\qp \}) > 0$ that was made in \cite{AJ-tap}.) In particular, one can choose $\eps=\eps_N\downarrow 0$ slowly enough such that 
\begin{equation*}
\lim_{N\to\infty}
\e\bigl\la \bigl( F_N^{\tap}(m(\sigma)) - F_N \bigr)^2\bigr\ra = 0.
\end{equation*}
One can combine this with the construction in \cite{JAU} to define a collection of pure states satisfying $F_N^{\tap}(m)\approx F_N$ in addition to all the properties obtained in \cite{JAU}. In other words, one can enforce the TAP free energy representation within pure states  in addition to all the usual properties of pure states. 

\medskip
\noindent
\textbf{Acknowledgements.} W.-K. C. thanks Antonio Auffinger, David Belius, and Nicola Kistler for fruitful discussions during the early stage of this work. Both authors thank the anonymous referees for the careful reading and providing several suggestions regarding the presentation of the paper. The research of W.-K. C. is partially supported by NSF DMS-1642207 and Hong Kong Research Grants Council GRF-14302515. D. P. is partially supported by NSERC.


\section{Zero-temperature Parisi formula for soft spins}\label{sec3}

Let $S$ be an arbitrary compact subset of $\mathbb{R}$ containing at least two elements. Recall the Gaussian process $X_N$ from \eqref{hamx}. We think of $X_N$ as a Hamiltonian defined on the product space $S^N$ with $\sum_{p\geq 2}\beta_p^2d(S)^p<\infty$ for $d(S):=2\max\{|s|:s\in S\}.$ Let $J$ be a continuous function on $S$. Consider the mixed $p$-spin model with a generalized external field $J$,
\begin{align*}
H_N^J(m)&=X_N(m)+h\sum_{i=1}^Nm_i+\sum_{i=1}^NJ(m_i),\,\,m\in S^N.
\end{align*}
This section is devoted to deriving the Parisi formula for the maximal energy of $H_N^J$ with self-overlap constraint. It will mainly be used in the proof of Theorem \ref{thm2}, but is also of independent interest. Our main result is stated in the following subsection.

\subsection{The statement of the formula}\label{sec2.1}

Let $\mathcal{D}$ be the convex hull of $\{s^2:s\in S\}$. For $u\in\mathcal{D}$, denote by $\mathcal{N}_u$ the space of all nonnegative right-continuous and nondecreasing functions on $[0,u)$ with finite $L^1(dx)$-norm. Define a functional $\mathcal{P}_u$ on $\mathbb{R}\times \mathcal{N}_u$ by
\begin{align}\label{pfinfinity}
\mathcal{P}_u(\lambda,\gamma)=\Phi_{u,\gamma}(0,h,\lambda)-\lambda u-\frac{1}{2}\int_0^u\xi''(s)s\gamma(s) ds,
\end{align}
where, for a given $\lambda,$ $\Phi_{u,\gamma}(s,x, \lambda)$ is defined as the weak solution of the following PDE,
\begin{align}
\label{pde2A}
\partial_s \Phi_{u,\gamma}&=-\frac{\xi''(s)}{2}\bigl(\partial_{xx}\Phi_{u,\gamma}+\gamma(s)\bigl(\partial_x\Phi_{u,\gamma}\bigr)^2\bigr),\,\,(s,x)\in [0,u)\times\mathbb{R}
\end{align}
with the boundary condition 
\begin{align}
\label{add:eq1}
\Phi_{u,\gamma}(u,x,\lambda)=f(x,\lambda):=\max_{m\in S}\bigl(m x+\lambda m^2+J(m) \bigr).
\end{align}
We show that the maximal energy of the Hamiltonian $H_N^J$ can be expressed as follows.

\begin{theorem} \label{thm1} For any $u\in \mathcal{D},$ we have that
	\begin{align}
	\label{thm1:eq1}
	\lim_{\varepsilon\downarrow 0}\lim_{N\rightarrow\infty} \e \max_{m\in S^N:|R(m,m)-u|<\varepsilon}\frac{H_N^J(m)}{N}=\inf_{(\lambda,\gamma)\in\mathbb{R}\times\mathcal{N}_u}\mathcal{P}_u(\lambda,\gamma).
	\end{align}
\end{theorem}
A few remarks are in order.

\begin{remark}\label{rmk1}\rm
	Theorem \ref{thm1} does not claim the existence of the minimizer. In Proposition \ref{add:lem4} below, a sufficient condition for the existence of the minimizer in terms of the directional derivative of the Parisi functional is established. 
\end{remark}

\begin{remark}\label{rmk2}\rm
Recall the Hamiltonian $H_N$ defined in \eqref{eq3}. If $S=\{-1,1\}$ and $J\equiv 0$, then $H_N=H_N^J$ and $\mathcal{D}=\{1\}.$ In this case, the above formula gives the Parisi formula for the maximum energy of $H_N,$ which agrees with the one previously obtained in Auffinger-Chen \cite{AC16}. 
\end{remark}

\begin{remark}\label{rmk3}\rm
When $S$ is a finite set, Theorem 1.2 of Jagannath-Sen \cite{JS17} established the Parisi formula for the maximum energy of the mixed even $p$-spin model (when $\xi$ is defined with $\beta_p=0$ for all odd $p\geq 3$) with general prior spin distributions \cite{Panchenko2005} by utilizing a generalization of the Parisi formula at positive temperature  \cite{Panchenko2015,Panchenk2015}. More precisely, their result states that when $S$ is finite,
	\begin{align*}
\lim_{\varepsilon\downarrow 0}\lim_{N\rightarrow\infty} \e \max_{m\in S^N:|R(m,m)-u|<\varepsilon}\frac{X_N(m)}{N}=\inf_{(\lambda,\nu)\in\mathbb{R}\times\mathcal{A}_u}\mathcal{P}_u'(\lambda,\nu).
\end{align*}
Here $\mathcal{A}_u$ is the space of all pairs $(\lambda,\nu)$ with $\lambda\in \mathbb{R}$ and $\nu(dt)=\gamma(t)dt+c\delta_u$ a measure defined on $[0,1]$ for $\gamma\in \mathcal{N}_u$ and $c\geq 0$. The functional $\mathcal{P}_u'(\lambda,\nu)$ has a form similar to \eqref{pfinfinity},
\begin{align*}
\mathcal{P}_u'(\lambda,\nu)=\Phi_{u,\nu}(0,0,\lambda)-\lambda u-\frac{1}{2}\int_0^u\xi''(s)s\nu(ds),
\end{align*}
where
\begin{align*}
\partial_s \Phi_{u,\nu}&=-\frac{\xi''(s)}{2}\bigl(\partial_{xx}\Phi_{u,\nu}+\gamma(s)\bigl(\partial_x\Phi_{u,\nu}\bigr)^2\bigr),\,\,(s,x)\in [0,u)\times\mathbb{R}
\end{align*}
and the boundary condition 
\begin{align*}
\Phi_{u,\nu}(u,x,\lambda)=\max_{m\in S}\Bigl(m x+\Bigl(\lambda +\frac{\xi''(u)}{2}c\Bigr)m^2 \Bigr).
\end{align*}
Since $\lambda$ varies over $\mathbb{R}$, if we make the change of variables $\lambda' =\lambda+\xi''(u)c/2$, the boundary condition becomes $\Phi_{u,\nu}(u,x,\lambda')=\max_{m\in S}(mx+\lambda'm^2)$ and the functional $\mathcal{P}_u'(\lambda,\nu)$ becomes
\begin{align*}
\mathcal{P}_u'(\lambda',\nu)=\Phi_{u,\nu}(0,0,\lambda')-\lambda' u-\frac{1}{2}\int_0^u\xi''(s)s\gamma(s)ds.
\end{align*}
These imply that 
$
\mathcal{P}_u'(\lambda,\nu)=\mathcal{P}_u(\lambda',\gamma),
$
so the formula of Theorem 1.2 in \cite{JS17} matches Theorem \ref{thm1} with the choice $J\equiv 0$ and $h=0.$ Thus, Theorem \ref{thm1} is a generalization of Theorem 1.2 in \cite{JS17} to not necessarily finite $S.$
\end{remark}

\begin{remark}\label{rmk4}\rm
	Although Theorem 1.2 in \cite{JS17} does not include the external field $J(m),$ adding $J(m)$ does not affect their proof in any way. Another difference is that Theorem 1.2 in \cite{JS17} was stated only for even-spin models, because the results in \cite{Panchenko2005, Panchenko2015} are proved only for even-spin models. However, this assumption is used in \cite{Panchenko2015} only in the Guerra-Talagrand upper bound. For general vector spin models, it is currently not known how to extend this upper bound to models that include odd $p$-spin terms. However, when the spins are one dimensional (as is the case we consider here, $S\subseteq \Reals$), Talagrand's positivity principle can be used to prove the upper bound for general mixed $p$-spin models in  exactly the same way as for the classical case of $\pm 1$ spins (see e.g. Sections 3.3 and 3.4 in \cite{Pan}).
\end{remark}

\subsection{Proof of the Parisi formula}

Our proof of Theorem \ref{thm1} as well as the argument of Theorem \ref{ThHTtapLT} rely on a very useful representation of the Parisi PDE in terms of the stochastic optimal control problem. We first gather some facts about the Parisi PDE and crucial properties of such representations. 

Let $0\leq a<b\leq 1.$ Suppose that $\phi$ is a real-valued Lipschitz function on $\mathbb{R}^2$ and $A$ is a nonnegative nondecreasing function on $[a,b)$ with right continuity and $\int_a^bA(s)ds<\infty.$ Consider the weak solution to the following PDE,
\begin{align}
\label{pde}
\partial_s\Psi(s,x,\lambda)&=-\frac{\xi''(s)}{2}\bigl(\partial_{xx}\Psi(s,x,\lambda)+A(s)\bigl(\partial_x\Psi(s,x,\lambda)\bigr)^2\bigr)
\end{align}
for $(s,x)\in [a,b)\times\mathbb{R}$ with the boundary condition $\Psi(b,x,\lambda)=\phi(x,\lambda).$ Note that the existence and uniqueness of the week solution to this PDE has appeared in many recent works \cite{CHL16,JS17,JT16} under many different boundary functions with Lipschitz property. An identical argument extends to the current setting, so the existence and uniqueness of the weak solution to \eqref{pde} are valid. Furthermore, it can be shown (see \cite{CHL16}) that $\Psi$ is differentiable with respect to $x$ up to any order on $[a,b)\times \mathbb{R}$ and these partial derivatives are uniformly bounded on $[a,b')\times \mathbb{R}$ for any finite $a<b'<b$. 

A key feature of the PDE \eqref{pde} is that it is a special case of the Hamilton-Jacobi-Bellman equation, induced by a linear problem of diffusion control \cite{FS}. In this case, $\Psi$ is known to possess a stochastic optimal control representation \cite{AC14,JS17}. More precisely, let $\mathcal{V}$ be the collection of all progressively measurable processes $v$ on $[a,b]$ with respect to the filtration generated by a standard Brownian motion $W$ and
\begin{align*}
\p\bigl(\max_{a\leq s\leq b}|v(s)|\leq 1\bigr)=1.
\end{align*}
We equip $\mathcal{V}$ with the metric $d(v,v'):=\bigl(\e\int_a^b(v(s)-v'(s))^2ds\bigr)^{1/2}.$ The PDE solution
 $\Psi$ can be written as a stochastic optimal control representation,
	\begin{align}
	\begin{split}\label{add:lem6:eq1}
	\Psi(a,x,\lambda)&=\max_{v\in \mathcal{V}}\e\Bigl[\phi\Bigl(x+\int_a^b\xi''(s)A(s)v(s)ds+\int_a^b\sqrt{\xi''(s)}dW(s),\lambda\Bigr)\\
	&\qquad\qquad\qquad\qquad-\frac{1}{2}\int_a^b\xi''(s)A(s)v(s)^2ds\Bigr],
	\end{split}
	\end{align}
	where the maximum of the right-hand side is attained by $$v(s)=\partial_x\Psi(s,X(s),\lambda),$$ where $X(s)$ is the solution of 
	$$
	dX(s)=\xi''(s)A(s)\partial_x\Psi(s,X(s),\lambda)ds+\sqrt{\xi''(s)}dW(s),\,\,X(a)=x.
	$$
	In addition, we also have that for any $a< s<b,$
	\begin{align}
	\begin{split}
	\label{add:eq2}
	d\partial_x\Psi(s,X(s),\lambda)&=\sqrt{\xi''(s)}\partial_{x}^2\Psi(s,X(s),\lambda)dW(s),\\
	d\partial_{x}^2\Psi(s,X(s),\lambda)&=-\xi''(s)A(s)\bigl(\partial_{x}^2\Psi(s,X(s),\lambda)\bigr)^2ds+\sqrt{\xi''(s)}\partial_{x}^3\Psi(s,X(s),\lambda)dW(s).
		\end{split}
	\end{align}
    The verification of \eqref{add:lem6:eq1} and \eqref{add:eq2} follows by a standard application of Ito's formula, which can be found, e.g., in \cite{AC14,AC16}. We do not reproduce the proof here.
    
    \smallskip
    
     We now turn to the proof of Theorem \ref{thm1}. The case when $S$ is finite was established in \cite{JS17}. We consider only the case that $S$ is an infinite set. We divide our discussion into two cases: $u$ in the interior $\mbox{Int}(\mathcal{D})$ of $\mathcal{D}$  and $u$ on the boundary $\mbox{Bd}(\mathcal{D})$ of $\mathcal{D}$. 

{\noindent \bf Proof of Theorem \ref{thm1} assuming $u\in \mbox{Int}(\mathcal{D})$.} Let $u\in Int(\mathcal{D})$ be fixed.
	For any $\delta>0,$ let $S_\delta$ be a finite subset of $S$ such that the $\delta$-neighborhood of $S_\delta$ covers $S$ and $|S_\delta|=O(\delta^{-1})$. Denote by $\mathcal{D}_\delta$ the convex hull of $\{s^2:s\in S_\delta\}.$ In addition, let us fix two distinct $s_1,s_2\in S$ such that $s_1^2< u< s_2^2$ and let us include them into $S_\delta$ for all $\delta.$ Then $u\in Int(\mathcal{D}_\delta)$ and  
	\begin{equation}\label{bee}
	d(u,\mathcal{D}_\delta^c)\geq b(u):=\min(s_2^2-u,u-s_1^2)>0
	\end{equation}
	for all $\delta>0.$ For notational convenience, set 
	$$A_{N,a}=\{m\in S^N:|R(m,m)-u|<a\},\,\,a>0.$$ 
	For any $m\in A_{N,\delta}$, the construction of $S_\delta$ ensures that there exists a $\pi(m)\in S_\delta^N$ such that 
	\begin{align*}
	\|m-\pi(m)\|<K\delta,
	\end{align*}
	where $\|\,\cdot\,\|$ is the Euclidean distance scaled by $1/\sqrt{N}$, and
	\begin{align*}
	\bigl|R\bigl(\pi(m),\pi(m)\bigr)-u\bigr|<K\delta,
	\end{align*}
	where $K$ is a fixed constant independent of $\delta$ and $N.$
    Our proof is completed by the following upper and lower bounds:
	
	{\noindent \bf Upper bound.} Define an auxiliary free energy
	\begin{align*}
	F_{\varepsilon,N}^\delta(\beta)&=\frac{1}{\beta N}\e\log \int_{m\in S_\delta^N:|R(m,m)-u|\leq \varepsilon}\exp\beta H_N^J(m)\mu^{\otimes N}(dm),\,\,\forall\beta,\varepsilon>0,
	\end{align*}
	where $\mu$ is a uniform probability measure on $S_\delta.$ The Parisi functional associated to this free energy is defined as follows. Let
	\begin{align*}
	f_\beta^\delta(x,\lambda):=\frac{1}{\beta}\log\int_{S_\delta}\exp \beta\bigl(xm+\lambda m^2+J(m)\bigr)d\mu.
	\end{align*}
	Let $\mathcal{N}_{\beta,u}$ be the collection of $\gamma\in \mathcal{N}_u$ with $\gamma(u)\leq \beta.$ Define the functional $\mathcal{P}_{\beta,u}^\delta$ on $\mathbb{R}\times\mathcal{N}_{\beta,u}$ by
	\begin{align*}
	\mathcal{P}_{\beta,u}^\delta(\lambda,\alpha)=\Phi_{\beta,u,\gamma}^\delta(0,h,\lambda)-\lambda u-\frac{1}{2}\int_0^u\xi''(s)s\gamma(s)\,ds,
	\end{align*}
	where $\Phi_{\beta,u,\gamma}^\delta$ is defined through the weak solution to the following PDE,
	\begin{align*}
	\partial_s\Phi_{\beta,u,\gamma}^\delta&=-\frac{\xi''(s)}{2}\bigl(\partial_{xx}\Phi_{\beta,u,\gamma}^\delta+\gamma(s) \bigl(\Phi_{\beta,u,\gamma}^\delta\bigr)^2\bigr),\,\,(s,x)\in [0,u)\times \mathbb{R},
	\end{align*}
	with the boundary condition $\Phi_{\beta,u,\gamma}^\delta(u,x,\lambda)=f_\beta(x,\lambda)$. Note that by applying \eqref{add:lem6:eq1} to both $\Phi_{u,\gamma}$ and $\Phi_{\beta,u,\gamma}^\delta$, we can immediately deduce that $\mathcal{P}_{\beta,u}^\delta\leq \mathcal{P}_u$ since $f_\beta^\delta\leq f.$ Using Dudley's entropy integral bound (see, e.g. \cite{BLM}), there exists a constant $C$ such that for any $N\geq 1$ and $\delta>0,$ 
	\begin{align}\label{dud}
	\e \max_{m,m'\in S^N:\|m-m'\|<\delta}\bigl|H_N^J(m)-H_N^J(m')\bigr|\leq C\delta N.
	\end{align}
	From \eqref{dud}, an application of the triangle inequality implies that
	\begin{align}
	\begin{split}
	\label{lem1:proof:eq3}
	\e \max_{m\in A_{N,\delta}}\frac{H_N^J(m)}{N}&\leq \e \max_{m\in A_{N,\delta}}\frac{\bigl|H_N^J(m)-H_N^J(\pi(m))\bigr|}{N}+\e \max_{m\in A_{N,\delta}}\frac{H_N^J(\pi(m))}{N}\\
	&\leq CK\delta+\e \max_{m\in A_{N,\delta}}\frac{H_N^J(\pi(m))}{N}\\
	&\leq CK\delta+\e \max_{m\in A_{N,K\delta}\cap S_\delta^N}\frac{H_N^J(m)}{N}\\
	&\leq CK\delta+F_{K\delta,N}^\delta(\beta)+\frac{\log |S_\delta|}{\beta}.
	\end{split}
	\end{align}
	To control the second term $F_{K\delta,N}^\delta(\beta)$, we use the standard Guerra replica symmetry breaking bound \cite{Guerra} (see, e.g., the proof of \cite[Lemma 2]{Panchenko2015}), which implies that, for any $N\geq 1$,
	\begin{align*}
	F_{K\delta,N}^\delta(\beta)&\leq \beta L\delta+\lambda \delta+\mathcal{P}_{\beta,u}^\delta(\lambda,\gamma)
	\end{align*}
	for $(\lambda,\gamma)\in\mathbb{R}\times\mathcal{N}_{\beta,u}.$
	Here $L$ is a constant depending only on $\xi'(1).$ Thus,
	\begin{align*}
	F_{K\delta,N}^\delta(\beta)&\leq  \beta L\delta+\lambda \delta+\mathcal{P}_u(\lambda,\gamma)
	\end{align*}
	for $(\lambda,\gamma)\in \mathbb{R}\times\mathcal{N}_{\beta,u}.$ Now if we take $\beta=1/\sqrt{\delta}$ and recall that $|S_\delta|=O(\delta^{-1})$, this inequality together with \eqref{lem1:proof:eq3} leads to
	\begin{align*}
	\lim_{\delta\downarrow 0}\limsup_{N\rightarrow\infty}\e \max_{m\in A_{N,\delta}}\frac{H_N^J(m)}{N}\leq \inf_{\mathbb{R}\times\mathcal{N}_u}\mathcal{P}_u(\lambda,\gamma).
	\end{align*}
	
	{\noindent \bf Lower bound.} Let $\mathcal{P}_u^\delta$ be defined in the same way as $\mathcal{P}_u$ with $S=S_\delta$, so that $f^{\delta}(x,\lambda):=\max_{m\in S_\delta}\bigl(m x+\lambda m^2+J(m) \bigr).$ Obviously,
	\begin{align*}
	\lim_{\varepsilon\downarrow 0}\liminf_{N\rightarrow\infty}\e \max_{m\in A_{N,\varepsilon}}\frac{H_N^J(m)}{N}&\geq \lim_{\varepsilon\downarrow 0}\liminf_{N\rightarrow\infty}\e \max_{m\in A_{N,\varepsilon}\cap S_\delta}\frac{H_N^J(m)}{N}=\inf_{\mathbb{R}\times\mathcal{N}_u}\mathcal{P}_u^\delta(\lambda,\gamma),
	\end{align*}
	where the last equality was proved in Theorem 1.2 in \cite{JS17} (see Remarks \ref{rmk3} and \ref{rmk4}). To finish the proof, we show that
		$$
		\lim_{\delta\downarrow 0}\inf_{\mathbb{R}\times\mathcal{N}_u}\mathcal{P}_u^\delta(\lambda,\gamma)\geq\inf_{\mathbb{R}\times\mathcal{N}_u}\mathcal{P}_u(\lambda,\gamma).
		$$
		First, recall from \cite{JS17} (proof of Theorem 3.9 and Lemma 3.10) that, for each $\delta>0,$ the infimum of $\mathcal{P}_u^\delta$ over $\mathbb{R}\times\mathcal{N}_u$ agrees with the infimum of $\mathcal{P}_u^\delta$ over all $(\lambda,\gamma)\in\mathbb{R}\times\mathcal{N}_u$  satisfying 
		\begin{align}
		\label{add:lem1:eq4}
		|\lambda|&\leq\frac{C}{d(u,\mathcal{D}_\delta^c)}\leq \frac{C}{b(u)}\,\,\mbox{and}\,\,\int_0^u\gamma(s)ds\leq C,
		\end{align}
		for some constant $C$ depending only on $\xi'(1)$ and $\xi''(1)$, and where we also used (\ref{bee}). Observe that there exists a constant $K>0$ such that
		\begin{align}\label{add:lem1:eq6}
		0\leq f(x,\lambda)-f^{\delta}(x,\lambda)&\leq K\delta \bigl(|x|+|\lambda|+1\bigr)+w(\delta)
		\end{align}
		for all $\delta>0$ and $\lambda,x\in \mathbb{R}$, where 
		$$w(\delta):=\sup_{m,m'\in S:|m-m'|<\delta}|J(m)-J(m')|.$$ 
        An application of the representation \eqref{add:lem6:eq1} together with \eqref{add:lem1:eq4} and \eqref{add:lem1:eq6} yields
	\begin{align*}
	0\leq \mathcal{P}_u(\lambda,\gamma)-\mathcal{P}_u^\delta(\lambda,\gamma)\leq K\delta\Bigl(C\xi''(1)+\e \Bigl|\int_0^u\sqrt{\xi''(s)}dW(s)\Bigr|+Cb(u)^{-1}+1\Bigr)+w(\delta)
	\end{align*}
	whenever $(\lambda,\gamma)$ satisfies \eqref{add:lem1:eq4}. Letting $\delta\downarrow 0$ completes our proof.
	\qed
	
	\medskip
	
	{\noindent \bf Proof of Theorem \ref{thm1} assuming $u\in \mbox{Bd}(\mathcal{D})$.} Suppose that the boundary of $\cal D$ is $\mbox{Bd}(\mathcal{D})=\{a,b\}$ for some $0\leq a<b.$ Take $u\in \{a,b\}$ and let $S_0=\{s\in S: s^2=u\}$. If we denote $s_0=\sqrt{u}$ then either $S_0=\{s_0\}$ or $S_0=\{-s_0,+s_0\}.$ Since $u$ is on the boundary, the constraint $|R(m,m)-u|<\eps$ in (\ref{thm1:eq1}) implies that the number of coordinates $m_i$ such that $|m_i^2-u|\geq \sqrt{\eps}$ is at most $\sqrt{\eps}N.$ This means that such $m$ are close to vectors in $S_0^N$ as $\eps\downarrow 0$ with respect to $\|\,\cdot\,\|$, which is the Euclidean norm normalized by $\sqrt{N}$. As in (\ref{dud}), this implies that
$$
\lim_{\varepsilon\downarrow 0}\lim_{N\rightarrow\infty} \e \max_{m\in S^N:|R(m,m)-u|<\varepsilon}\frac{H_N^J(m)}{N}	
=
\lim_{N\rightarrow\infty} \e \max_{m\in S_0^N}\frac{H_N^J(m)}{N}.
$$	
By Theorem \ref{thm1} for finite sets (Theorem 1.2 in \cite{JS17}), we know that this limit equals 
$$
\lim_{N\rightarrow\infty} \e \max_{m\in S_0^N}\frac{H_N^J(m)}{N}
=
\inf_{(\lambda,\gamma)\in \mathbb{R}\times \mathcal{N}_u}\mathcal{P}_u^0(\lambda,\gamma),
$$
where $\mathcal{P}_u^0(\lambda,\gamma)$ is defined as in (\ref{pfinfinity}) with the boundary condition $\max_{m\in S_0}\bigl(m x+\lambda m^2+J(m) \bigr).$ Since shifting the boundary condition by a constant results in shifting the PDE solution in (\ref{pde2A}) by the same constant, we can move the term $-\lambda u$ in (\ref{pfinfinity}) into the boundary conditions and rewrite	
\begin{align*}
\mathcal{P}^0_u(\lambda,\gamma)=\Phi_{u,\gamma}(0,h,\lambda)-\frac{1}{2}\int_0^u\xi''(s)s\gamma(s) ds,
\end{align*}
where $\Phi_{u,\gamma}(0,h, \lambda)$ satisfies the PDE (\ref{pde2A}) with the boundary condition 
\begin{align*}
f_0(x,\lambda)=\max_{m\in S_0}\bigl(m x+\lambda (m^2-u)+J(m) \bigr).
\end{align*}
Notice that $m^2=u$ for $m\in S_0$, so $f_0(x,\lambda):=\max_{m\in S_0}\bigl(m x+J(m) \bigr)$ and $\mathcal{P}^0_u(\lambda,\gamma)$ does not depend on $\lambda.$ As a result,
$$
\lim_{N\rightarrow\infty} \e \max_{m\in S_0^N}\frac{H_N^J(m)}{N}
=
\inf_{(\lambda,\gamma)\in \mathbb{R}\times \mathcal{N}_u}\mathcal{P}_u^0(\lambda,\gamma)
=
\inf_{\gamma\in \mathcal{N}_u}\mathcal{P}_u^0(0,\gamma).
$$
To finish the proof of the theorem, it remains to show that
	\begin{align}\label{add:eq8}
	\inf_{\gamma\in \mathcal{N}_u}\mathcal{P}_u^0(0,\gamma)
	=
	\inf_{\mathbb{R}\times \mathcal{N}_u}\mathcal{P}_u(\lambda,\gamma).
	\end{align}
Observe that, by shifting $-\lambda u$ as above, we can rewrite $\mathcal{P}_u(\lambda,\gamma)$  in (\ref{pfinfinity})	as
\begin{align*}
\mathcal{P}_u(\lambda,\gamma)=\Phi_{u,\gamma}(0,h,\lambda)-\frac{1}{2}\int_0^u\xi''(s)s\gamma(s) ds,
\end{align*}
where $\Phi_{u,\gamma}(0,h, \lambda)$ is the solution of the PDE (\ref{pde2A}) with the boundary condition 
\begin{align*}
f(x,\lambda)=\max_{m\in S}\bigl(m x+\lambda (m^2-u)+J(m) \bigr).
\end{align*}
Since $f(x,\lambda) \geq f_0(x,\lambda),$ by monotonicity of the solution of the PDE (\ref{pde2A}) with respect to the boundary condition,
$$
	\inf_{ \mathcal{N}_u}\mathcal{P}_u^0(0,\gamma)
	=
	\inf_{\mathbb{R}\times \mathcal{N}_u}\mathcal{P}_u^0(\lambda,\gamma)
	\leq
	\inf_{\mathbb{R}\times \mathcal{N}_u}\mathcal{P}_u(\lambda,\gamma).
$$
In the other direction, since $u$ is on the boundary, $m^2-u$ is either positive or negative for all $m\in S$. Suppose for certainty that $m^2- u\leq 0$. In this case, we write
$$
\inf_{(\lambda,\gamma)\in\mathbb{R}\times \mathcal{N}_u}\mathcal{P}_u(\lambda,\gamma)
\leq
\inf_{(\lambda,\gamma)\in\mathbb{R}^+\times \mathcal{N}_u}\mathcal{P}_u(\lambda,\gamma)
\leq
\inf_{\gamma\in\mathcal{N}_u} \lim_{\lambda\to+\infty}\mathcal{P}_u(\lambda,\gamma).
$$ 
For $\lambda>0$ and $m^2- u\leq 0$, the boundary condition can be written as
\begin{align*}
f(x,\lambda)=\max_{m\in S}\bigl(m x - \lambda |m^2-u|+J(m) \bigr)
\end{align*}
and it is clear that, for any $x\in \Reals$,
$$
\lim_{\lambda\to+\infty} f(x,\lambda) = \max_{m\in S_0}\bigl(m x +J(m) \bigr) = f_0(x,0).
$$
On the other hand,
\begin{align*}
\max_{m\in S_0}\bigl(m x +J(m) \bigr)\leq f(x,\lambda)\leq \max_{m\in S}\bigl(m x +J(m) \bigr),
\end{align*}
so the function $f(x,\lambda)$ is bounded uniformly in $\lambda$ by const$\,\times(1+|x|).$ Therefore, the variational representations \eqref{add:lem6:eq1} and the bounded convergence theorem imply that
	\begin{align*}
	\lim_{\lambda\to+\infty}\mathcal{P}_u(\lambda,\gamma)=\mathcal{P}_u^0(0,\gamma)
	\end{align*}
	and, therefore,
	$\inf_{\mathbb{R}\times \mathcal{N}_u}\mathcal{P}_u(\lambda,\gamma)\leq \inf_{\mathcal{N}_u}\mathcal{P}_u^0(0,\gamma).$
	This finishes the proof.
\qed

\section{TAP free energy representation}\label{sec4}

The aim of this section is to establish Theorem \ref{thm2}. In Subsection \ref{sub4.1} we express the maximum of the TAP free energy with self-overlap constraint via a Parisi-type formula using Theorem \ref{thm1}. Next, in Subsection \ref{sub3.2} we derive the directional derivative of the Parisi functional and extract useful criteria for the optimizers. The proof of Theorem \ref{thm2} is presented in Subsection \ref{sub4.3}.

\subsection{Parisi formula for the restricted TAP free energy}\label{sub4.1}

Recall the Hamilton $H_N$ from \eqref{eq3}, the rate function $I$ from \eqref{eq-7}, the correction term $C$ from \eqref{eq-8}, and the TAP free energy $F_N^{\tap}$ from \eqref{tap}. For $u\in [0,1]$, the maximum TAP free energy with self-overlap constraint is defined as
\begin{align*}
F^{\tap}_{\varepsilon,N}(u)&=\max_{m\in [-1,1]^N:|R(m,m)-u|<\varepsilon}F_N^{\tap}(m).
\end{align*} 
Throughout the remainder of the paper, we set $S=[-1,1]$ and $J(x)=-I(x)$. Applying Theorem~\ref{thm1} with these choices shows that, for $u\in(0,1),$
\begin{align*}
\lim_{\varepsilon\downarrow 0}\lim_{N\rightarrow\infty}\e \max_{m\in [-1,1]^N:|R(m,m)-u|<\varepsilon}\Bigl(\frac{H_N(m)}{N}-\frac{\sum_{i=1}^NI(m_i)}{N}\Bigr)&=\inf_{\mathbb{R}\times\mathcal{N}_u}\mathcal{P}_u(\lambda,\gamma),
\end{align*}
where $\mathcal{P}_u$ is defined in \eqref{pfinfinity}. Consequently, the constrained TAP free energy can be written as a Parisi formula, for $u\in[0,1],$
\begin{align}\label{pfTAP}
F^{\tap}(u):=\lim_{\varepsilon\downarrow 0}\lim_{N\rightarrow\infty} \e F^{\tap}_{\varepsilon,N}(u)&=\inf_{\mathbb{R}\times\mathcal{N}_u}\mathcal{P}_u(\lambda,\gamma)+C(u).
\end{align}

\subsection{Directional derivative}\label{sub3.2}

We now compute the direction derivative of the Parisi functional $\mathcal{P}_u$ following the framework of Chen \cite{C14}. Let  $(\lambda_0,\gamma_0)$, $(\lambda,\gamma)\in \mathbb{R}\times\mathcal{N}_u$. For any $\theta\in [0,1]$, set
\begin{align*}
\lambda_\theta&=(1-\theta)\lambda_0+\theta\lambda,
\gamma_\theta=(1-\theta)\gamma_0+\theta\gamma.
\end{align*}

\begin{proposition}
	\label{add:lem2}
	Let $(\lambda_0,\gamma_0),(\lambda,\gamma)\in\mathbb{R}\times\mathcal{N}_u$ with $|\lambda_0|<1/4$ and $\gamma_0(u-),\gamma(u-)<\infty$. Then
	\begin{align*}
	\frac{d}{d\theta}\mathcal{P}_u(\lambda_\theta,\gamma_\theta)\Bigl|_{\theta=0}&=\frac{1}{2}\int_0^u\xi''(s)(\gamma(s)-\gamma_0(s))\bigl(\e v_0^2(s)-s\bigr)ds+(\lambda-\lambda_0)\bigl(\e v_0^2(u)-u\bigr),
	\end{align*}
	where the derivative is from the right-hand side of $0$ and $v_0(s):=\partial_x\Phi_{u,\gamma_0}(s,X_0(s))$ for $0\leq s\leq u$, where $X_0$ is the solution of
	\begin{align*}
	dX_0(s)&=\xi''(s)\gamma(s) \partial_x\Phi_{u,\gamma_0}(s,X_0(s))ds+{\xi''(s)}^{1/2}dW(s),\,\,X_0(0)=h.
	\end{align*}
	\end{proposition}
	
	The following proposition gives a criterion for reaching the optimality of $\mathcal{P}_u$.
	
\begin{proposition}
	\label{add:lem4}
	If there exists some $(\lambda_0,\gamma_0)\in \mathbb{R}\times\mathcal{N}_u$ with $|\lambda_0|<1/4$ and $\gamma_0(u-)<\infty$ such that
	\begin{align}\label{add:lem4:eq1}
	\frac{d}{d\theta}\mathcal{P}_u(\lambda_\theta,\gamma_\theta)\Bigl|_{\theta=0}\geq 0
	\end{align}
	for all $(\lambda,\gamma)\in \mathbb{R}\times\mathcal{N}_u$ with $\gamma(u-)<\infty$, then $(\lambda_0,\gamma_0)$ minimizes $\mathcal{P}_u$ over $\mathbb{R}\times\mathcal{N}_u.$
\end{proposition}

In the remainder of this section we establish Propositions \ref{add:lem3} and \ref{add:lem4}. Recall from Subsection~\ref{sub4.1} that $S=[-1,1]$ and $J=-I$. The boundary condition \eqref{add:eq1} is given by
\begin{align*}
f(x,\lambda)&=\max_{m\in [-1,1]}\bigl(mx+\lambda m^2-I(m)\bigr).
\end{align*}

\begin{lemma}\label{add:lem3}
For any $(x,\lambda)\in \mathbb{R}\times (-1/4,1/4),$
	\begin{align}\label{add:lem3:eq1}
	f(x,\lambda)&=\log 2\mathrm{ch}\bigl(x+2\lambda m(x,\lambda)\bigr)-\lambda m(x,\lambda)^2,
	\end{align}
	where $m(x,\lambda)$ is the unique solution to 
		\begin{align}\label{add:lem3.1:eq1}
		m=\myth(x+2\lambda m)
		\end{align}
		for $m\in[-1,1].$ In addition, $f$ is twice continuously differentiable on $\mathbb{R}\times (-1/4,1/4)$ with
	\begin{align}
	\begin{split}\label{add:lem3:eq2}
	\partial_xf(x,\lambda)&=m(x,\lambda),\\
	\partial_\lambda f(x,\lambda)&=m(x,\lambda)^2
	\end{split}
	\end{align}
	and uniformly bounded second derivatives.
\end{lemma}
{\noindent \bf Proof.}
	Set $
	g(x,\lambda,m)=mx+\lambda m^2-I(m).
	$
	A direct computation yields that 
	\begin{align*}
	\frac{\partial}{\partial m}g(x,\lambda,m)&=x+2\lambda m-\frac{1}{2}\log \frac{1+m}{1-m}=0
	\end{align*}
	if and only if $m$ satisfies \eqref{add:lem3.1:eq1}. The equation \eqref{add:lem3:eq1} can be checked by applying \eqref{add:lem3.1:eq1}. The fact that $m(x,\lambda)$ is the unique solution to \eqref{add:lem3.1:eq1} is because the map $m\in[-1,1]\mapsto\myth(x+2\lambda m)$ is a contraction for $|\lambda|<1/4.$ The rest of the statements of Lemma \ref{add:lem3} can be easily derived by using the equation \eqref{add:lem3.1:eq1} and performing implicit differentiation.
\qed

\medskip
{\noindent \bf Proof of Proposition \ref{add:lem2}.} Note that the boundary condition $f$ of $\Phi_{u,\gamma_\theta}$ is Lipschitz on $\mathbb{R}^2$. From \eqref{add:lem6:eq1},
	\begin{align}\label{add:lem2:proof:eq1}
	\Phi_{u,\gamma_\theta}(0,h,\lambda_\theta)&=\max_{v\in \mathcal{V}}\e \Bigl[f\Bigl(h+\int_0^u\xi''\gamma_\theta vds+\int_0^u{\xi''}^{1/2}dW,\lambda_\theta\Bigr)-\frac{1}{2}\int_0^u\xi''\gamma_\theta v^2ds\Bigr],
	\end{align}
	where the maximum is achieved by the process $v_\theta(s):=\partial_x\Phi_{u,\gamma_\theta}(s,X_\theta(s),\lambda_\theta),$ where
	\begin{align}\label{add:eq-3}
	dX_\theta&=\xi''\gamma_\theta \partial_x\Phi_{u,\gamma_\theta}(s,X_\theta(s),\lambda_\theta)ds+\sqrt{\xi''}dW,\,\,X_\theta(0)=h.
	\end{align}	
	From these, 
	\begin{equation}\label{Puthgamma}
	\mathcal{P}_u(\lambda_\theta,\gamma_\theta)=\max_{v\in \mathcal{V}} g(v,\theta)=g(v_\theta,\theta),
	\end{equation}
	where
	\begin{align*}
	g(v,\theta)&:=\e \Bigl[f\Bigl(h+\int_0^u\xi''\gamma_\theta vds+\int_0^u\sqrt{\xi''}dW,\lambda_\theta\Bigl)-\frac{1}{2}\int_0^u\xi''\gamma_\theta (v^2+s)ds\Bigr].
	\end{align*}
	By assumption, $|\lambda_0|<1/4$, so we can choose $\theta'>0$ small enough so that $|\lambda_\theta|<1/4$ for all $\theta\in [0,\theta'].$ For such $\theta,$ the function $f(x,\lambda_\theta)$ satisfies all the properties in Lemma \ref{add:lem3} and, in particular, it is twice continuously differentiable with bounded derivatives on $\mathbb{R}\times(-1/4,1/4).$ Using the finiteness of $\gamma_0(u-)$ and $\gamma(u-)$, a standard argument (see, e.g., \cite{JT16}) yields that $\partial_x^i\Phi_{u,\gamma_\theta}(s,x,\lambda_\theta)$ for $i=1,2$ is continuous and uniformly bounded. Thus, using \eqref{add:eq-3} and Gronwall's inequality one can easily check that $\theta\mapsto v_\theta$ is continuous. Finally, by adapting \cite[Lemma 2]{C14}, $\mathcal{P}_u(\lambda_\theta,\gamma_\theta)$ in (\ref{Puthgamma}) is right-differentiable at zero and its derivative is equal to $\partial_\theta g(v_\theta,\theta)$. It remains to compute this derivative. A direct computation using Lemma \ref{add:lem3} yields, for $|\lambda_\theta|<1/4$,
	\begin{align*}
	\begin{split}
	\partial_\theta g(v,\theta)&=\e \Bigl[m\Bigl(h+\int_0^u\xi''\gamma_\theta v\,ds+\int_0^u\sqrt{\xi''}dW,\lambda_\theta\Bigr)\int_0^u\xi''(\gamma-\gamma_0)v\,ds\Bigr]\\
	&+\e \Bigl[m\Bigl(h+\int_0^u\xi''\gamma_\theta v\,ds+\int_0^u\sqrt{\xi''}dW,\lambda_\theta\Bigr)^2\Bigr](\lambda-\lambda_0)\\
	&-\frac{1}{2}\int_0^u\xi''(\gamma-\gamma_0) (\e v^2+s)ds.
	\end{split}
	\end{align*}
	If we plug in $v=v_\theta$ then, by (\ref{add:eq-3}),
	\begin{align*}
	\begin{split}
	\partial_\theta g(v_\theta,\theta)=&\,\, \e \Bigl[m(X_\theta(u),\lambda_\theta)\int_0^u\xi''(\gamma-\gamma_0)v_\theta\,ds\Bigr]\\
	&\,+\,\e m(X_\theta(u),\lambda_\theta)^2 (\lambda-\lambda_0)-\frac{1}{2}\int_0^u\xi''(\gamma-\gamma_0) (\e v_\theta^2+s)ds.
	\end{split}
	\end{align*}
	Since the boundary condition gives $\Phi_{u,\gamma_\theta}(u,x,\lambda_\theta) = f(x, \lambda_\theta),$ by (\ref{add:lem3:eq2}) we get that $\partial_x\Phi_{u,\gamma_\theta}(u,x,\lambda_\theta) = m(x,\lambda_\theta)$ and, therefore, $m(X_\theta(u),\lambda_\theta) = \partial_x\Phi_{u,\gamma_\theta}(u,X_\theta(u),\lambda_\theta) = v_\theta(u)$ and
	\begin{align}
	\begin{split}\label{add:lem3:proof:eq1}
	\partial_\theta g(v_\theta,\theta)=&\,\, \e v_\theta(u) \int_0^u\xi''(\gamma-\gamma_0)v_\theta\,ds 
	\\
	&\,+\,\e v_\theta(u)^2 (\lambda-\lambda_0)-\frac{1}{2}\int_0^u\xi''(\gamma-\gamma_0) (\e v_\theta^2+s)ds.
	\end{split}
	\end{align}
	From the first equation of \eqref{add:eq2},
	\begin{align*}
	dv_\theta(s)=\xi''(s)^{1/2}\partial_{xx}\Phi_{u,\gamma_\theta}(s,X_\theta(s),\lambda_\theta)\,dW(s),
	\end{align*}
	which implies that $\e v_\theta(u)v_\theta(s) = \e v_\theta(s)^2$ and this allows us to rewrite \eqref{add:lem3:proof:eq1} as
	\begin{align*}
	\partial_\theta g(v_\theta,\theta)
	&=\int_0^u\xi''(\gamma-\gamma_0)\e v_\theta^2 ds+\e v_\theta(u)^2(\lambda-\lambda_0)-\frac{1}{2}\int_0^u\xi''(\gamma-\gamma_0) (\e v_\theta^2+s)ds\\
	&=\frac{1}{2}\int_0^u\xi''(\gamma-\gamma_0)\bigl(\e v_\theta^2-s\bigr) ds+\e v_\theta(u)^2(\lambda-\lambda_0).
	\end{align*} 
	 This finishes the proof.
\qed

\medskip
{\noindent \bf Proof of Proposition \ref{add:lem4}.}
	Observe that the boundary condition $f$ defined in \eqref{add:eq1} is convex in $\mathbb{R}^2.$ By the virtue of the representation \eqref{add:lem2:proof:eq1} as well as the definition of $\mathcal{P}_u$ in Theorem \ref{thm1}, one can easily see that $\mathcal{P}_u$ is a convex functional on $\mathbb{R}\times\mathcal{N}_u.$ As a consequence, the condition \eqref{add:lem4:eq1} yields that $\mathcal{P}_u(\lambda,\gamma)\geq \mathcal{P}_u(\lambda_0,\gamma_0)$ for all $(\lambda,\gamma)\in \mathbb{R}\times\mathcal{N}_u$ with $\gamma(u-)<\infty$. The restriction $\gamma(u-)<\infty$ can then be removed by an approximation argument.
\qed

\smallskip

\subsection{Proof of Theorem \ref{thm2}}\label{sub4.3}

Recall the Parisi functional $\mathcal{P}$ from \eqref{pf}. Also recall that the Parisi measure $\ap(ds)$ is the unique minimizer of $\mathcal{P}$ and $\qp$ is the largest point in the support of $\ap(ds).$ We divide our proof into two parts: $u\in(\qp,1]$ and $u=\qp.$ The first case will establish $F^{\tap}(u)\leq F$ for all $u\in (\qp,1]$, while the second case aims to obtain $F^{\tap}(\qp)=F.$ These together finishes our proof of Theorem \ref{thm2}.

\subsubsection{Case I: $u\in (\qp,1]$}

Obviously, from \eqref{add:lem3:eq1}, $f(x,0)=\log\mathrm{ch} (x).$ From \eqref{pfTAP} and the comment below it,
\begin{align}
\begin{split}
\label{add:eq3}
F^{\tap}(u)&\leq \inf_{\gamma\in \mathcal{N}_u}\mathcal{P}_u(0,\gamma)+C(u)\\
&\leq \inf_{\gamma\in \mathcal{M}_u}\mathcal{P}_u(0,\gamma)+C(u),
\end{split}
\end{align}
where we recall that $\mathcal{M}_u$ is the space of all cumulative distribution function on $[0,u].$
Observe that if $\alpha:=\gamma 1_{[0,u)}+1_{[u,1]}$ for $\gamma\in \mathcal{M}_u,$ then
\begin{align*}
\Phi_\alpha(0,h)&=\Phi_{u,\gamma}(0,h,0)+\frac{1}{2}\bigl(\xi'(1)-\xi'(u)\bigr)
\end{align*}
and
\begin{align*}
\int_0^1\xi''(s)\alpha(s)sds&=\int_0^u\xi''(s)\gamma(s)sds+\int_u^1\xi''(s)sds\\
&=\int_0^u\xi''(s)\gamma(s)sds+\xi'(1)-\xi'(u)u-\bigl(\xi(1)-\xi(u)\bigr).
\end{align*}
Thus, combining these equations implies
$
\mathcal{P}_u(0,\gamma)=\mathcal{P}(\alpha)-C(u)
$
and, by \eqref{add:eq3},
$
F^{\tap}(u)\leq \mathcal{P}(\alpha)
$
for all $\alpha$ of the form $\gamma 1_{[0,u)}+1_{[u,1]}$ for $\gamma\in \mathcal{M}_u.$ Since the minimum of $\mathcal{P}$ over $\mathcal{M}$ is attained by ${\ap}$ and $u>\qp$, this and the above inequality imply that $F^{\tap}(u)\leq F$ for all $u\in(\qp,1].$

\subsubsection{Case II: $u=\qp$}\label{5.2}

The second case is more involved. It is based on the tools we developed in Subsection \ref{sub3.2}. First, we recall that analogous results of Propositions \ref{add:lem2} and \ref{add:lem4} are also valid for the functional $\mathcal{P}$ defined in \eqref{pf}.

\begin{proposition}
	\label{prop1} Let $\alpha_0\in\mathcal{M}$ be fixed.
	For any $\alpha\in\mathcal{M}$, define $\alpha_\theta=(1-\theta)\alpha_0+\theta\alpha$ for $\theta\in[0,1]$. The directional derivative of $\mathcal{P}$ can be computed as
	\begin{align*}
	\frac{d}{d\theta}\mathcal{P}(\alpha_\theta)\Big|_{\theta=0}=\frac{1 }{2}\int_0^1\xi''(s)(\e w_{0}(s)^2-s)(\alpha(s)-\alpha_0(s))ds
	\end{align*}
	for all $\alpha\in\mathcal{M},$ where the derivative is from the right-hand side of $0$. Here $w_0(s)=\partial_x\Phi_{\alpha_0}(s,Y_0(s))$, where $Y_0$ is the solution of
	$$dY_0=\xi''\alpha_0\partial_x\Phi_{\alpha_0}(s,Y_0)ds+\sqrt{\xi''} dW, Y_0(0)=h.$$ In addition, $\alpha_0$ is the Parisi measure if and only if $\frac{d}{d\theta}\mathcal{P}(\alpha_\theta)\big|_{\theta=0}\geq 0$ for all $\alpha\in\mathcal{M}.$
\end{proposition}

Next, we need a key lemma that will be used to establish the equality between the restricted TAP free energy and the free energy. 
\begin{lemma}
	\label{prep:lem1}
	Let $\alpha_0\in \mathcal{M}$ be fixed. Assume that the support $\Omega$ of $\alpha_0(ds)$ is contained in $[0,1)$. Let $g$ be a continuous function on $[0,1)$ such that, whenever $0\in \Omega$, $g(0)=0$. If 
	\begin{align*}
	\int_0^1g(s)(\alpha(s)-\alpha_0(s))ds\geq 0
	\end{align*}
	for all $\alpha\in\mathcal{M}$ then, for any $u\in \Omega,$
	\begin{align}\label{prep:lem1:eq2}
	\int_0^{u}g(s)(\gamma(s)-\alpha_0(s))ds\geq 0
	\end{align}
	for all $\gamma\in \mathcal{N}_u.$
\end{lemma}

{\noindent \bf Proof.}
First we claim that for any $r,r_1,r_2\in \Omega,$ 
\begin{align}
\label{prelim:eq1}
g(r)&=0
\end{align}
and
\begin{align}
\label{prelim:eq2}
\int_{r_1}^{r_2}g(s)ds=0.
\end{align}
Define 
$
\bar{g}(r)=\int_r^1 g(s)ds.
$
By Fubini's theorem,
\begin{align*}
\int_0^1(\alpha(s)-\alpha_0(s))g(s)ds=\int_0^1\Bigl(\int_0^s(\alpha-\alpha_0)(dr)\Bigr)g(s)ds=\int_0^1\bar{g}(r)(\alpha-\alpha_0)(dr).
\end{align*}
Thus, the given assumption implies
\begin{align*}
\int_0^1\bar g(r)\alpha_0(dr)\leq \int_0^1\bar g(r)\alpha(dr).
\end{align*}
Since this holds for all $\alpha\in\mathcal{M}$, it is equivalent to $\alpha_0(A)=1$ for
$$
A:=\Bigl\{s\in[0,1]\,:\,\bar g(s)=\min_{r\in[0,1]} \bar g(r) \Bigr\}.
$$
This implies that $A$ contains $\Omega$ and, as a result, 
$$\int_{r_1}^{r_2}g(s)ds=\bar g(r_1)-\bar g(r_2)=0.$$ 
If $r>0$ then, from the optimality of $\bar g$, $g(r)=\bar g'(r)=0.$ If $r=0,$ then $g(r)=0$ by the given assumption. These finish the proof of our claim. 

Now we turn to the proof of \eqref{prep:lem1:eq2}. Let $u\in \Omega$ be fixed. Denote $u'=\min\, \Omega$. Since $(u',u)\setminus \Omega$ is open, this set can be written as the countable union of disjoint open intervals $(I_n)_{n\in\mathbb{N}}.$ On each such interval, since $\alpha_0$ is a constant, \eqref{prelim:eq2} implies that
\begin{align*}
\int_{I_n}g(s)\alpha_0(s)ds=0.
\end{align*}
This and \eqref{prelim:eq1} yield
\begin{align}\label{prelim:eq3}
\int_0^ug(s)\alpha_0(s)ds&=\int_0^{u'}g(s)\alpha_0(s)ds+\sum_{n\in \mathbb{N}}\int_{I_n}g(s)\alpha_0(s)ds=0,
\end{align}
where the first integral vanishes because $\alpha_0=0$ on $[0,u').$ Now, consider $\gamma\in \mathcal{N}_u$ with $\gamma(u)<\infty.$ Let $c>\gamma(u)/\alpha_0(u).$ Set
\begin{align*}
\alpha=\frac{\gamma}{c}1_{[0,u)}+\alpha_01_{[u,1]}\in\mathcal{M}.
\end{align*}  
From the given assumption and \eqref{prelim:eq3},
\begin{align*}
0&\leq \int_0^1 g(s)(\alpha(s)-\alpha_0(s))ds=\int_0^ug(s)(\alpha(s)-\alpha_0(s))ds=\int_0^ug(s)\alpha(s)ds=\frac{1}{c}\int_0^ug(s)\gamma(s)ds.
\end{align*}
Another use of \eqref{prelim:eq3} gives
\begin{align*}
0&\leq \int_0^ug(s)\gamma(s)ds=\int_0^ug(s)\gamma(s)ds-\int_0^u g(s)\alpha_0(s)ds=\int_0^ug(s)(\gamma(s)-\alpha_0(s))ds.
\end{align*}
This yields \eqref{prep:lem1:eq2} in the case $\gamma(u)<\infty$. This condition can be removed by an obvious approximation argument.
\qed

\medskip
The proof of the second case proceeds as follows. Let $u=\qp.$ Denote $\alpha_0=\ap$ on $[0,1]$ and $\gamma_0=\ap$ on $[0,u].$ Also, let $\lambda_0=0.$ Note that 
\begin{align}\label{add:eq5}
\log 2+\Phi_{\alpha_0}(s,x)&=\Phi_{u,\gamma_0}(s,x,0)+\frac{1}{2}\bigl(\xi'(1)-\xi'(u)\bigr)
\end{align}
for $(s,x)\in[0,u]\times\mathbb{R}.$ This implies that
\begin{align*}
\partial_x\Phi_{\alpha_0}(s,x)&=\partial_x\Phi_{u,\gamma_0}(s,x,0).
\end{align*}
Recall $v_0$ and $w_0$ from Propositions \ref{add:lem2} and \ref{prop1}. From this equation, 
\begin{align}
\label{add:eq4}
v_0=w_0\,\,\mbox{on $[0,u]$}.
\end{align}
From Proposition \ref{prop1},
$$
\int_0^1g(s)(\alpha(s)-\alpha_0(s))ds\geq 0,\,\,\forall \alpha\in \mathcal{M},
$$
where \begin{align*}
g(s):=\xi''(s)\bigl(\e w_0(s)^2-s\bigr)
\end{align*} 
for $s\in [0,1].$ Note that if $h\neq 0$ then $\Omega$ does not contain $0$ (see \cite[Theorem 14.12.1]{Talbook2}); if $h=0$ then $0\in \Omega$ (see \cite[Theorem 1]{AC131}) and in this case since $\partial_x\Phi_{\alpha_0}(0,\,\cdot\,)$ is an odd function, we see $\e w_0(0)^2=\partial_x\Phi_{\alpha_0}(0,0)^2=0$. Thus, $g(0)=0$ whenever $0\in \Omega.$ Therefore, Lemma \ref{prep:lem1} implies that
$$
\int_0^ug(s)(\gamma(s)-\gamma_0(s))ds\geq 0
$$   
for all $\gamma\in \mathcal{N}_u$. Also, we see that $g(u)=0$ from \eqref{prelim:eq1}. If $u\neq 0$, then this equation forces $\e w_0(u)^2=u$. If $u=0,$ it must be true that $h=0$  and  $\e w_0(0)^2=0.$ From these, \eqref{add:eq4}, and Proposition~\ref{add:lem2}, we obtain 
\begin{align*}
\frac{d}{d\theta}\mathcal{P}_u(\lambda_\theta,\gamma_\theta)\Bigl|_{\theta=0}&=\int_0^ug(s)(\gamma(s)-\gamma_0(s))ds+(\lambda-\gamma_0)(\e v_0^2(u)-u)\\
&=\int_0^ug(s)(\gamma(s)-\gamma_0(s))ds \geq 0
\end{align*} 
for all $\gamma\in \mathcal{N}_u$. Thus, by Proposition \ref{add:lem4}, $(\lambda_0,\gamma_0)$ is a minimizer of $\mathcal{P}_u$ over $\mathbb{R}\times\mathcal{N}_u$ and a direct computation using the definitions of $\alpha_0,\gamma_0$ and the equation \eqref{add:eq5} gives $\mathcal{P}_u(\lambda_0,\gamma_0)=\mathcal{P}(\alpha_0)-C(u).$ This means that
$$
F^{\tap}(u)=\mathcal{P}_u(\lambda_0,\gamma_0)+C(u)=\mathcal{P}(\alpha_0)=F.
$$ 
This finishes our proof.


\section{TAP free energy at high temperature}

In this section, we will prove Theorems \ref{ThHTtap} and \ref{ThHTent}. In the proof of Theorem \ref{ThHTtap}, we will use Theorem 1.7.1 from the book of Talagrand \cite{Talbook1} that was proved for the pure SK model corresponding to $\xi(s)=\beta^2 s^2/2$ and which was stated for $\beta<1/2.$ However, the proof of Theorem 1.7.1 is based on a simple cavity computation and can be easily extended to general $\xi$ under the condition (\ref{HTcond}), with one difference that the explicit rate of convergence $O(1/N)$ in \cite{Talbook1} would now depend on the rate of convergence in (\ref{HTcond}). This is because the condition (\ref{HTcond}) automatically holds for the cavity Gibbs measures, by Gronwall's lemma. For simplicity, below we will stick to the case of $\xi(s)=\beta^2 s^2/2$.

\medskip
\noindent
\textbf{Proof of Theorem \ref{ThHTtap}.} Note that $C(u)=\beta^2(1-u)^2/4$. From the condition \eqref{HTcond}, the correction term can be easily handled by using the mean value theorem and the Cauchy-Schwarz inequality,
\begin{align*}
\lim_{N\rightarrow\infty}\e \bigr(C(\qe)-C(q)\bigl)^2&=\lim_{N\rightarrow\infty}\e \bigr(C(\la R_{1,2}\ra)-C(q)\bigl)^2\\
&\leq \frac{\beta^2}{2}\lim_{N\rightarrow\infty}\e \bigl(\la R_{1,2}\ra-q\bigr)^2\\
&=\frac{\beta^2}{2}\lim_{N\rightarrow\infty}\e\bigl\la \bigl(R_{1,2}-q\bigr)^2\bigr\ra=0.
\end{align*}
To handle the entropy term
$$
I_N(m) = \frac{1}{N}\sum_{i\leq N} I(\la\sigma_i\ra),
$$ 
recall from Theorem 1.7.1 in \cite{Talbook1} that $\la\sigma_i\ra$ can be approximated in distribution by $\myth(\beta z_i \sqrt{q}+h)$ for standard Gaussian $z_i \sim N(0,1)$ and, moreover, $\la\sigma_i\ra$ and $\la\sigma_j\ra$ are asymptotically uncorrelated for any $i\not = j.$ (By symmetry between sites, the distribution of these pairs are the same for all $(i,j).$) Therefore, $\lim_{N\to\infty}\mathrm{Var}\,(I_N(m)) = 0,$ so the term $I_N(m)$ concentrates and its expected value is approximated by $\e I(\myth(\beta z \sqrt{q}+h)).$ Using the identity
$$
I(\myth(x))=\frac{1+\myth(x)}{2}\log\frac{(1+\myth(x))}{2}+\frac{1-\myth(x)}{2}\log\frac{(1-\myth(x))}{2}
= x\,\myth(x)-\log 2\ch(x)
$$
and integration by parts, we can rewrite
\begin{align*}
\e I(\myth(\beta z \sqrt{q}+h))
&= \e (\beta z \sqrt{q}+h)\myth(\beta z \sqrt{q}+h)-\e \log2 \mathrm{ch}(\beta z \sqrt{q}+h)
\\
& = h\,\e \myth(\beta z \sqrt{q}+h) +\beta^2 q\e(1- \myth^2(\beta z \sqrt{q}+h)) - \e  \log 2\mathrm{ch}(\beta z \sqrt{q}+h)
\\
& = h\,\e \myth(\beta z \sqrt{q}+h) +\beta^2 q(1-q) - \e  \log2 \mathrm{ch}(\beta z \sqrt{q}+h),
\end{align*}
where in the last step we used that $\e\myth^2(\beta z \sqrt{q}+h) = q.$

Next, let us denote $\beta_0 = \beta/\sqrt{2}$ and consider
$$
\frac{H_N(m)}{N} = 
\frac{\beta_0}{N^{3/2}}\sum_{i,j} g_{ij} \la \sigma_i\ra \la\sigma_j\ra
+ \frac{h}{N}\sum_{i\leq N} \la\sigma_i\ra.
$$
By the same reasoning as above, the second term concentrates around
$$
\frac{h}{N}\sum_{i\leq N} \la\sigma_i\ra \approx h\,\e \myth(\beta z \sqrt{q}+h).
$$
To show that the first term concentrates, let us compute its first and second moments. Let us denote
$$
H_N^{\mathrm{SK}}(\sigma^1,\sigma^2) = \frac{1}{\sqrt{N}}\sum_{i,j} g_{ij}\sigma_i^1\sigma_j^2
$$
and observe that
$$
\frac{1}{N} \e H_N^{\mathrm{SK}}(\sigma^1,\sigma^2) H_N^{\mathrm{SK}}(\sigma^3,\sigma^4)
=
R_{1,3} R_{2,4}.
$$
Then the usual Gaussian integration by parts gives that the first moment is
\begin{align*}
\e \frac{1}{N^{3/2}} \sum_{i,j} g_{ij} \la \sigma_i\ra \la\sigma_j\ra = \frac{1}{N}\e \la H_N^{\mathrm{SK}}(\sigma^1,\sigma^2) \ra
=
2 \beta_0\, \e\la R_{1,2} - R_{1,3} R_{2,3}\ra
\approx 2\beta_0 q(1-q),
\end{align*}
where in the last step we used (\ref{HTcond}). To compute the second moment, we rewrite it as
\begin{align*}
&\frac{1}{N^2}\e \la H_N^{\mathrm{SK}}(\sigma^1,\sigma^2) H_N^{\mathrm{SK}}(\sigma^3,\sigma^4) \ra\\
& =
\frac{\beta_0}{N} \e\la H_N^{\mathrm{SK}}(\sigma^3,\sigma^4)
(2R_{1,2} +R_{1,3} R_{2,3}+ R_{1,4} R_{2,4} - 4R_{1,5} R_{2,5}) \ra
+ \frac{1}{N}\e\la R_{1,3} R_{2,4}\ra\\
& \approx 2\beta_0 q(1-q) \frac{1}{N} \e\la H_N^{\mathrm{SK}}(\sigma^3,\sigma^4) \ra\\
&\approx (2\beta_0 q(1-q))^2,
\end{align*}
where in the last two steps we used the high temperature condition (\ref{HTcond}) and the fact that we already computed  the first moment in the previous step. All together, this shows that
$$
\frac{H_N(m)}{N} 
\approx \beta^2 q(1-q) + h\,\e \myth(\beta z \sqrt{q}+h).
$$
Finally, it remains to recall the formula for the free energy from \cite[Proposition 1.6.8]{Talbook1} under the condition (\ref{HTcond}),
$$
F_N\approx \frac{\beta^2}{4}(1-q)^2 + \e\log 2\ch(\beta z \sqrt{q}+h).
$$
Combining all the terms finishes the proof.
\qed

\medskip
\noindent
\textbf{Proof of Theorem \ref{ThHTent}.} 
It is well known that both free energy $F_N$ and average energy $\la H_N(\sigma)\ra/N$ in (\ref{FGibbsRep}) concentrate, so the entropy $\la \log G_N(\sigma)\ra/N$ also concentrates. Since
$$
\frac{\la H_N(\sigma) \ra}{N} = 
\frac{\beta_0}{N} \la X_N(\sigma)\ra+ \frac{h}{N}\sum_{i\leq N} \la\sigma_i\ra,
$$
where we again use the notation $\beta_0=\beta/\sqrt{2},$ and
$$
\frac{\beta_0}{N} \e\la X_N(\sigma)\ra = \beta_0^2(1-\e\la R_{1,2}^2\ra) \approx \frac{\beta_0^2}{2}(1-q^2),
$$
we get from (\ref{FGibbsRep}) and the fact that the Gibbs measure $G_N$ reaches the optimality of \eqref{FGibbsRep} that
\begin{align*}
\frac{1}{N}\e \la \log G_N(\sigma)\ra
&\approx
\frac{\beta^2}{2}(1-q^2)+h\e \mathrm{th}(\beta z \sqrt{q}+h)-\frac{\beta^2}{4}(1-q)^2-\e\log 2\mathrm{ch}(\beta z\sqrt{q}+h).
\end{align*}
To finish the proof it remains to compare this with $I_N(\la\sigma\ra)$ in the proof of Theorem \ref{ThHTtap} above.
\qed

\section{TAP free energy for pure states at low temperature}

Throughout this section, we assume that the model is generic. We start with a simple technical lemma that will be useful below.

\begin{lemma}\label{CondD}
	For any $\eps>0,$
	$$
	\lim_{\delta\downarrow 0}\liminf_{N\to\infty} \,\e G_N(\sigma \,:\, W_N(\sigma)\geq \delta) = 1.
	$$
\end{lemma}
\textbf{Proof.}
First of all, by decreasing $\eps$, we decrease the set $\Sigma_N(\sigma)$ and its measure $W_N(\sigma).$ Therefore, we can assume that $\qp -\eps$ is the point of continuity of the Parisi measure $\ap ,$ which is the limiting distribution of the overlap $R_{1,2}.$ As a result,
\begin{align*}
\lim_{N\to\infty} \,\e \la W_N(\sigma)^k \ra 
&= \lim_{N\to\infty} \,\e \la I(R(\sigma,\rho^\ell)\geq \qp -\eps, 1\leq \ell\leq k) \ra 
\\
&= \e \la I(\bsigma\cdot \brho^\ell \geq \qp -\eps, 1\leq \ell\leq k) \ra 
= \e \la W(\bsigma)^k \ra,
\end{align*}
where on the right hand side the average $\la\,\cdot\,\ra$ is with respect to the asymptotic Gibbs measure $G$, and $W(\bsigma) = G(\Sigma(\bsigma))$, where 
$$
\Sigma(\bsigma) =\bigl\{\brho \in L^2([0,1]) \,:\, \bsigma\cdot \brho\geq \qp -\eps\bigr\}.
$$
Moreover, since $G$ lives on the sphere of radius $\sqrt{\qp },$ we can restrict $\brho$ to $\{\|\brho\|^2=\qp \}.$ For generic mixed $p$-spin model, the asymptotic Gibbs measure is described by the Ruelle probability cascades (see \cite{Pan}) and, in particular, we can rewrite
$$
\e \la W(\bsigma)^k \ra = \e \sum_{n\geq 1} v_n^{k+1},
$$
where the sequence of weights $(v_n)_{n\geq 1}$ has the Poisson-Dirichlet distribution $PD(\zeta)$ with $\zeta=\p(\bsigma\cdot \brho>\qp -\eps) = 1-\ap ([0,\qp -\eps]).$ Since these weights have continuous distribution, we can approximate $I(x\geq \delta)$ by polynomials to conclude that
$$
\lim_{N\to\infty} \,\e \la I(W_N(\sigma)\geq \delta) \ra 
= \e \sum_{n\geq 1} v_n I(v_n\geq \delta).
$$
Letting $\delta\downarrow 0$ and using that $\sum_{n\geq 1} v_n=1$ finishes the proof.
\qed

\subsection{Computation of the first term}

Next, we will control the mixed $p$-spin term of the Hamiltonian, excluding the external field. Recall the Gaussian process $X_N$ from \eqref{hamx}. 

\begin{lemma}\label{LemFT}
	For generic $p$-spin models,
	\begin{equation*}
	\lim_{\eps\downarrow 0}\limsup_{N\to\infty}
	\e\Bigl\la \Bigl( \frac{X_N(\la\rho \ra_\sigma)}{N}
	-\e\Bigl\la \frac{X_N(\la\rho \ra_\sigma)}{N} \Bigr\ra\Bigr)^2\Bigr\ra = 0
	\end{equation*}
	and
	\begin{equation}\label{FTmom}
	\lim_{\eps\downarrow 0}\limsup_{N\to\infty}\Bigl|
	\e\Bigl\la \frac{X_N(\la\rho \ra_\sigma)}{N} \Bigr\ra
	- \Bigl(\xi'(\qp ) - \theta(\qp )- \int \xi(q) \,d\ap (q)\Bigr)
	\Bigr|
	=0
	\end{equation}
	for $\theta(s):=s\xi'(s)-\xi(s).$
\end{lemma}
\textbf{Proof.}
Express $X_N=\sum_{p\geq 2} \beta_p H_{N,p}.$ Our assumption $\sum_{p\geq 2} 2^p \beta_p^2 <\infty$ implies that $\sum_{p\geq 2} \beta_p <\infty,$ so it is enough to prove concentration for each $H_{N,p}$. We will do this for $p=2$ for simplicity of notation, because the general case differs only in the number of various indices involved. Set $$
A_\ell = \bigl\{R(\sigma,\sigma^\ell)\geq \qp -\eps\bigr\}.
$$ 
Let us start by computing the first moment
$$
\e\Bigl\la \frac{H_{N,2}(\la\rho \ra_\sigma)}{N} \Bigr\ra
=
\sum_{i,j}\frac{1}{N^{3/2}} \e\Bigl\la g_{ij} \sigma_i^1\sigma_j^2 I(A_1\cap A_2) \frac{1}{W_N(\sigma)^2}\Bigr\ra.
$$
Integration by parts (see Section A.4 in \cite{Talbook1}) gives
$$
\sum_{i,j}\frac{\beta_2}{N^{2}} \e\Bigl\la \sigma_i^1\sigma_j^2 I(A_1\cap A_2) \frac{1}{W_N(\sigma)^2}
\Bigl(
\sigma_i\sigma_j + \sigma_i^1\sigma_j^1 + \sigma_i^2\sigma_j^2 - \sigma^3_i\sigma_j^3 - 2\sigma^3_i\sigma_j^3 \frac{I(A_3)}{W_N(\sigma)}
\Bigr)
\Bigr\ra.
$$
Notice that, in the last two terms $\sigma^3_i\sigma_j^3$, one is unconstrained and the other one is constrained by the event $A_3= \{R(\sigma,\sigma^3)\geq \qp -\eps\}$, which comes from integrating by parts of $\sum_{\rho\in \Sigma_N(\sigma)} \exp H_N(\rho),$ which is the numerator of $W_N(\sigma).$ Summing over $i,j$, this equals
\begin{align*}
&\beta_2  \e\Bigl\la I(A_1\cap A_2) \frac{1}{W_N(\sigma)^2}
\Bigl(R(\sigma,\sigma^1)R(\sigma,\sigma^2) + 2R(\sigma^1,\sigma^2) \\
&\qquad- R(\sigma^1,\sigma^3)R(\sigma^2,\sigma^3) - 2R(\sigma^1,\sigma^3)R(\sigma^2,\sigma^3) \frac{I(A_3)}{W_N(\sigma)}
\Bigr)
\Bigr\ra.
\end{align*}
Let us fix $\delta>0$, denote $\phi = I(W_N(\sigma)\geq \delta)$, and break the above expectation into two terms -- with the factor $\phi$ and factor $1-\phi = I(W_N(\sigma)< \delta).$ We bound the second one by $6\beta_2 \e\la 1-\phi\ra$, which, by Lemma \ref{CondD}, satisfies
$$
\lim_{\delta\downarrow 0}\limsup_{N\to\infty} 6\beta_2 \e\la 1-\phi\ra = 0.
$$
The first term with the factor $\phi$ is
\begin{align*}
&\beta_2  \e\Bigl\la I(A_1\cap A_2) \frac{\phi}{W_N(\sigma)^2}
\Bigl(R(\sigma,\sigma^1)R(\sigma,\sigma^2) + 2R(\sigma^1,\sigma^2) \\
&\qquad- R(\sigma^1,\sigma^3)R(\sigma^2,\sigma^3) - 2R(\sigma^1,\sigma^3)R(\sigma^2,\sigma^3) \frac{I(A_3)}{W_N(\sigma)}
\Bigr)
\Bigr\ra.
\end{align*}
For simplicity, let us assume that $\eps>0$ is such that $\qp -\eps$ is the point of continuity of the Parisi measure $\ap $ (otherwise, one can use another approximation argument as $\eps\downarrow 0$). Because of the indicator $\phi = I(W_N(\sigma)\geq \delta),$ we can approximate $1/W_N(\sigma)^n$ by polynomials of $W_N(\sigma)$ on the interval $[\delta,1]$ and then, as in the proof of Lemma \ref{CondD}, use replicas to show that the above converges to the same expression (denoting $R(\bsigma^\ell,\bsigma^{\ell'})=\bsigma^\ell\cdot\bsigma^{\ell'}$),
\begin{align}
\begin{split}\label{add:eq-2}
&\beta_2  \e\Bigl\la I(A_1\cap A_2) \frac{\phi}{W(\bsigma)^2}
\Bigl(R(\bsigma,\bsigma^1)R(\bsigma,\bsigma^2) + 2R(\bsigma^1,\bsigma^2)\\
&\qquad - R(\bsigma^1,\bsigma^3)R(\bsigma^2,\bsigma^3) - 2R(\bsigma^1,\bsigma^3)R(\bsigma^2,\bsigma^3) \frac{I(A_3)}{W(\bsigma)}
\Bigr)
\Bigr\ra,
\end{split}
\end{align}
where everything is now expressed in terms of the asymptotic Gibbs measure $G.$ Using the constraints $A_\ell$, ultrametricity of the support of $G$ (see \cite{Pultra} or \cite{Pan}), and the fact that $\qp $ is the largest point in the support of the limiting distribution of the overlap, the above differs from
\begin{align}\label{add:eq-4}
\beta_2  \e\Bigl\la I(A_1\cap A_2) \frac{\phi}{W(\bsigma)^2}
\Bigl((\qp )^2 + 2\qp  - R(\bsigma,\bsigma^3)^2 - 2(\qp )^2 \frac{I(A_3)}{W(\bsigma)}
\Bigr)
\Bigr\ra
\end{align}
by at most $L\eps.$ In fact, all the overlaps in \eqref{add:eq-2} are in the interval $[\qp -\eps,\qp ]$, except for the overlaps $R(\bsigma^\ell,\bsigma^3)$ in the second to last term (where $\bsigma^3$ is not constrained), which, by constraint $A_\ell$ and ultrametricity, differ from $R(\bsigma,\bsigma^3)$ by at most $\eps.$ Now write the equation \eqref{add:eq-4} as
$$
\beta_2  \e\bigl\la \phi \bigl(2\qp  - (\qp )^2 - R(\bsigma,\bsigma^3)^2 \bigr)\bigr\ra
$$
and, writing $\phi=1-(1-\phi)$, this equals to
$$
\beta_2\bigl(2 \qp  - (\qp )^2 - \e \la R(\bsigma^1,\bsigma^2)^2 \ra \bigr)
=
\beta_2\bigl(2 \qp  - (\qp )^2 - \int q^2 \,d\ap (q) \bigr)
$$
up to an error $4\beta_2\e\la1-\phi\ra$, which goes to zero as $\delta\downarrow 0.$ Thus, we showed that
\begin{equation*}
\limsup_{N\to\infty}\Bigl|
\e\Bigl\la \frac{H_{N,2}(\la\rho \ra_\sigma)}{N} \Bigr\ra
- \beta_2\bigl(2 \qp  - (\qp )^2 - \int q^2 \,d\ap (q) \bigr)
\Bigr|
\leq L\eps.
\end{equation*}
The same computation for mixed $p$-spin Hamiltonian gives (\ref{FTmom}).

Similarly, we can compute the second moment. If we denote $A^j = \cap_{\ell\leq j}A_\ell$ then
$$
\e\Bigl\la \Bigl(\frac{H_{N,2}(\la\rho \ra_\sigma)}{N}\Bigr)^2\Bigr\ra
=
\sum_{i,j,k,\ell}\frac{1}{N^{3}} \e\Bigl\la g_{ij}g_{k\ell} \sigma_i^1\sigma_j^2 \sigma_k^3\sigma_\ell^4 I(A^4) \frac{1}{W_N(\sigma)^4}\Bigr\ra.
$$
This requires two integrations by parts. Let us fix $(i,j)\not = (k,\ell)$ (the terms $(k,\ell)=(i,j)$ will add up to lower order contribution) and integrate by parts with respect to $g_{ij}$ first to get $N^{-7/2}$ times
\begin{align*}
&
\beta_2 \e\Bigl\la g_{k\ell} \sigma_i^1\sigma_j^2 \sigma_k^3\sigma_\ell^4 I(A^4) \frac{1}{W_N(\sigma)^4}
\Bigl(
\sigma_i\sigma_j + \sum_{n\leq 4}\sigma_i^n\sigma_j^n -\sigma^5_i\sigma_j^5 
- 4\sigma^5_i\sigma_j^5 \frac{I(A_5)}{W_N(\sigma)}
\Bigr)
\Bigr\ra
\\
& =
\beta_2 \e\Bigl\la g_{k\ell} \sigma_i^1\sigma_j^2 \sigma_k^3\sigma_\ell^4 I(A^5) \frac{1}{W_N(\sigma)^5}
\Bigl(
\sigma_i\sigma_j + \sum_{n\leq 4}\sigma_i^n\sigma_j^n  
- 4\sigma^5_i\sigma_j^5
\Bigr)
\Bigr\ra
\\
&\hspace{1cm} -
\beta_2\e\Bigl\la g_{k\ell} \sigma_i^1\sigma_j^2  \sigma^5_i\sigma_j^5 \sigma_k^3\sigma_\ell^4 I(A^4) \frac{1}{W_N(\sigma)^4}
\Bigr\ra = I-II.
\end{align*}
Integrate by parts the first term $I$ with respect to $g_{k\ell}$ (the factor in front now is $N^{-4}$):
\begin{align*}
&
\beta_2^2\e\Bigl\la \sigma_i^1\sigma_j^2 \sigma_k^3\sigma_\ell^4 I(A^5) \frac{1}{W_N(\sigma)^5}
\Bigl(
\sigma_i\sigma_j + \sum_{n\leq 4}\sigma_i^n\sigma_j^n  
- 4\sigma^5_i\sigma_j^5
\Bigr)\times
\\
&\hspace{1cm} \times
\Bigl(
\sigma_k\sigma_\ell + \sum_{n\leq 5}\sigma_k^n\sigma_\ell^n -\sigma^6_k\sigma_\ell^6 
- 5\sigma^6_k\sigma_\ell^6 \frac{I(A_6)}{W_N(\sigma)}
\Bigr)
\Bigr\ra
\\
& =\,
\beta_2^2\e\Bigl\la \sigma_i^1\sigma_j^2 \sigma_k^3\sigma_\ell^4 I(A^6) \frac{1}{W_N(\sigma)^6}
\Bigl(
\sigma_i\sigma_j + \sum_{n\leq 4}\sigma_i^n\sigma_j^n  
- 4\sigma^5_i\sigma_j^5
\Bigr)
\Bigl(
\sigma_k\sigma_\ell + \sum_{n\leq 5}\sigma_k^n\sigma_\ell^n 
- 5\sigma^6_k\sigma_\ell^6
\Bigr)
\Bigr\ra
\\
&\hspace{1cm} -
\beta_2^2\e\Bigl\la \sigma_i^1\sigma_j^2 \sigma_k^3\sigma_\ell^4 \sigma^6_k\sigma_\ell^6 I(A^5) \frac{1}{W_N(\sigma)^5}
\Bigl(
\sigma_i\sigma_j + \sum_{n\leq 4}\sigma_i^n\sigma_i^n  
- 4\sigma^5_i\sigma_j^5
\Bigr)
\Bigr\ra.
\end{align*}
When we sum over $i,j,k,\ell,$ we can rewrite this in terms of overlaps. As in the argument for the first moment above, in the thermodynamic limit, constrained overlaps will be replaced by $\qp $ and in any unconstrained overlap one can replace a constrained replica by $\bsigma.$ One can check that, up to error of order $L\eps,$ we get
$$
\beta_2^2(2\qp -(\qp )^2)^2-\beta_2^2(2\qp -(\qp )^2) \int q^2 \,d\ap (q).
$$
When we integrate the second term, we get (up to factor $N^{-4}$)
\begin{align*}
&
\beta_2^2
\e\Bigl\la \sigma_i^1\sigma_j^2  \sigma^6_i\sigma_j^6 \sigma_k^3\sigma_\ell^4 I(A^4) \frac{1}{W_N(\sigma)^4}
\Bigl(
\sigma_k\sigma_\ell + \sum_{n\leq 5}\sigma_k^n\sigma_\ell^n -2\sigma^6_k\sigma_\ell^6 
- 4\sigma^6_k\sigma_\ell^6 \frac{I(A_6)}{W_N(\sigma)}
\Bigr)
\Bigr\ra.
\end{align*}
Again, one can check that in the limit this is, up to error terms of order $L\eps$,
$$
\beta_2^2 \int q^2 \,d\ap (q)
\Bigl(2\qp -(\qp )^2 - \e \int q^2 \,d\ap (q)\Bigr).
$$
Taking the difference of the two terms, we showed that
\begin{equation*}
\limsup_{N\to\infty}\Bigl|
\e\Bigl\la \Bigl(\frac{H_{N,2}(\la\rho \ra_\sigma)}{N}\Bigr)^2\Bigr\ra
- \beta_2^2\Bigl(2 \qp  - (\qp )^2 - \int q^2 \,d\ap (q)\Bigr)^2
\Bigr|
\leq L\eps.
\end{equation*}
This finishes the proof.
\qed

\subsection{Computation of the second term (and external field)}

In this section we will consider the combination of the external field and the second term in the TAP free energy. In other words, for $\varphi(x) := -hx + I(x),$ we will study
\begin{equation*}
Y_N(\sigma):=\frac{1}{N}\sum_{i\leq N}\varphi(m_i(\sigma)).
\end{equation*}
Denote by $\Phi(s,x)$ the solution of the Parisi PDE associated to the Parisi measure $\ap $. Recall from \eqref{add:lem6:eq1} that $\Phi(0,x)$ can be written as a stochastic optimal control problem and the optimum of this variational representation is reached by the process $\bigl(\partial_x\Phi(s,X(s))\bigr)_{0\leq s\leq 1}$ for $(X(s))_{s\in[0,1]}$ the solution to the following SDE,
\begin{align}\label{SDEField}
dX(s)&=\xi''(s)\ap (s)\partial_x\Phi(s,X(s))ds+\xi''(s)^{1/2}dW,\,\,X(0)=h.
\end{align}
We will prove the following.
\begin{lemma}\label{LemST}
	For generic $p$-spin models,
	\begin{equation}\label{LemSTconc}
	\lim_{\eps\downarrow 0}\limsup_{N\to\infty}
	\e\bigl\la \bigl( Y_N(\sigma)
	-\e\bigl\la Y_N(\sigma) \bigr\ra\bigr)^2\bigr\ra = 0
	\end{equation}
	and
	\begin{equation}\label{STmom}
	\lim_{\eps\downarrow 0}\limsup_{N\to\infty}\Bigl|
	\e\bigl\la Y_N(\sigma) \bigr\ra
	- \bigl(\e (X(\qp )-h)\mathrm{th}(X(\qp ))-\e \log 2\mathrm{ch}(X(\qp ))\bigr)
	\Bigr|
	=0.
	\end{equation}
\end{lemma}
\textbf{Proof.} 
The main tools of the proof will be the convergence of spin distributions from Chapter 4 in \cite{Pan} and an alternative representation using SDEs from \cite{AJ-tap}. We will approximate the function $\varphi(x) = -hx + I(x)$ by polynomials on $[-1,1]$, so we first consider the moments 
\begin{equation}\label{SecTYk}
M_k(\sigma):=\frac{1}{N}\sum_{i\leq N}m_i(\sigma)^k=\frac{1}{N}\sum_{i\leq N} \bigl\la\rho_i\bigr\ra_\sigma^k.
\end{equation}
Using replicas and symmetry between sites, we can rewrite 
$$
\e\la M_k(\sigma)\ra = 
\e\Bigl\la \sigma_1^1\cdots \sigma_1^k  \,I(R(\sigma,\sigma^\ell)\geq \qp -\eps, 1\leq \ell\leq k)\frac{1}{W_N(\sigma)^k}
\Bigr\ra.
$$
As in the proof of Lemma \ref{LemFT},  we fix $\delta>0$, denote $\phi = I(W_N(\sigma)\geq \delta)$, and break the above integral into two terms -- with the factor $\phi$ and factor $1-\phi = I(W_N(\sigma)< \delta).$ The second term is bounded by $\e\la 1-\phi\ra$, which, by Lemma \ref{CondD}, is small for small $\delta$.
The first term with the factor $\phi$ is
$$
\e\Bigl\la \sigma_1^1\cdots \sigma_1^k  \,I(W_N(\sigma)\geq \delta)\,I(R(\sigma,\sigma^\ell)\geq \qp -\eps, 1\leq \ell\leq k)\frac{1}{W_N(\sigma)^k}
\Bigr\ra.
$$
The factor $\phi$ allows us to approximate $1/W_N(\sigma)^{k}$ by polynomials of $W_N(\sigma)$, since the function $1/x^k$ is continuous on $[\delta,1],$ and then use replicas to express moments of $W_N(\sigma)$ in terms of functions of the overlaps, as in Lemma \ref{CondD}. We will again assume, for simplicity, that $\qp -\eps$ is the point of continuity of the Parisi measure $\ap $.

The results of Chapter 4 in \cite{Pan} were proved using regularizing perturbations developed in \cite{Pspins}; however, as we mentioned above, it was observed in \cite{AJ-tap} (Proposition A.7) that they hold for generic $p$-spin models without any perturbations. The results of Chapter 4 in \cite{Pan} allow us to compute rather explicitly the limits of joint moments of spins and overlaps in terms of the Parisi measure $\ap $. More precisely, the limit of such moments depends continuously on the asymptotic distribution of one overlap and, moreover, an explicit construction exists for distributions that concentrate of finitely many points.  As a result, one could compute these limits via discrete approximations $\ap '$ of the Parisi measure $\ap $. A different representation of this explicit construction was given in \cite{AJ-tap} that allows one to do some computations directly for arbitrary $\ap ,$ which we will utilize below. 

More precisely, the limit of the above expression is equal to
$$
\e\Bigl\la \myth(X(\bsigma^1))\cdots \,\myth(X(\bsigma^k))\, I(W(\bsigma)\geq \delta)
I(R(\bsigma,\bsigma^\ell)\geq \qp -\eps, 1\leq \ell\leq k)\frac{1}{W(\bsigma)^k}
\Bigr\ra,
$$
where the average $\la\,\cdot\,\ra$ is with respect to the Ruelle probability cascades $G$ corresponding to $\ap,$ and the stochastic process $X(\bsigma)$ indexed by the points $\bsigma$ in the support of $G$ was {explicitly described in} Theorem 4.4 in \cite{Pan} in the case when $\ap$ is discrete. We do not recall this explicit construction, since we will use an alternative description from \cite{AJ-tap} which applies to arbitrary $\ap$ and, moreover, we will need only one-point distribution of this process. If we denote by $\e_X$ the expectation with respect to $X$ given $G$ then, for fixed $\bsigma^1,\ldots, \bsigma^k$ satisfying the condition $\bsigma^\ell\cdot \bsigma\geq \qp -\eps,$, we will first use an approximation
$$
\bigl| \e_X \myth(X(\bsigma^1))\cdots \,\myth(X(\bsigma^k)) - \e_X \myth(X(\bsigma))^k \bigr| \leq k\sqrt{2\eps},
$$
which holds because one of the properties of the process $X$ states that $\e_X \myth(X(\bsigma^1))\myth(X(\bsigma^2)) = \bsigma^1\cdot\bsigma^2$ and, on the event $\bsigma^\ell\cdot \bsigma\geq \qp -\eps,$
$$
\bigl| \e_X \bigl(\myth(X(\bsigma^\ell)) - \myth(X(\bsigma)) \bigr)^2 \bigr| 
= \bsigma^\ell\cdot \bsigma^\ell + \bsigma\cdot \bsigma - 2\bsigma^\ell\cdot \bsigma \leq 2\eps.
$$
By Lemma 2.2 in \cite{AJ-tap}, the distribution of $X(\bsigma)$ for a fixed $\bsigma$ can be expressed via the SDE (\ref{SDEField}), so
$$
\e_X \myth(X(\bsigma))^k = \e \myth(X(\qp ))^k,
$$
where $X(\qp )$ is defined by (\ref{SDEField}). Thus, letting $\delta\downarrow 0$, we showed that
\begin{align}
\label{ad:eq-6}
\limsup_{N\to\infty}\bigl| \e\la M_k(\sigma)\ra - \e \myth(X(\qp ))^k \bigr| \leq k\sqrt{2\eps}.
\end{align}
The equation (\ref{STmom}) follows by approximation by polynomials and rewriting
$$
\varphi(\myth(X(\qp ))) = (X(\qp )-h)\mathrm{th}(X(\qp ))-\e \log 2\mathrm{ch}(X(\qp )).
$$
To show the concentration (\ref{LemSTconc}), one can similarly compute the second moment of $M_k(\sigma)$ (using symmetry between sites and ignoring diagonal terms)
\begin{align*}
\e\la M_k(\sigma)^2 \ra
= &\,\,
\e\Bigl\la \sigma_1^1\cdots \sigma_1^k 
I(R(\sigma,\sigma^\ell)\geq \qp -\eps,1\leq \ell\leq k)\frac{1}{W_N(\sigma)^k}\, \times
\\
&\hspace{7mm}
\times \sigma_2^{k+1}\cdots \sigma_2^{2k}
I(R(\sigma,\sigma^\ell)\geq \qp -\eps, k+1\leq \ell\leq 2k)\frac{1}{W_N(\sigma)^k}
\Bigr\ra + O(N^{-1}).
\end{align*}
The only difference is that, by the result in Chapter 4 in \cite{Pan}, the coordinates $\sigma_1^\ell$ are replaced by $\myth(X_1(\bsigma^\ell))$ and $\sigma_2^\ell$ are replaced by $\myth(X_2(\bsigma^\ell))$ for two independent copies $X_1$ and $X_2$ of the process $X(\bsigma)$. Otherwise, the same argument as above give
\begin{align}
\label{add:eq-7}
\limsup_{N\to\infty}\bigl| \e\la M_k(\sigma)^2\ra - \bigl(\e \myth(X(\qp ))^k \bigr)^2 \bigr| \leq 2k\sqrt{2\eps},
\end{align}
which implies the concentration (\ref{LemSTconc}).
\qed

\medskip

The following lemma relates the quantity (\ref{STmom}) to the limiting free energy $\PP(\ap )$ and the correction term $C(\qp)$. Denote
\begin{equation*}
\Delta =  \e (X(\qp )-h)\mathrm{th}(X(\qp ))-\e \log 2\mathrm{ch}(X(\qp )).
\end{equation*}

\begin{lemma}\label{lem6}
	We have
	\begin{equation*}
	\Delta =-C(\qp)-\PP(\ap ) + \int_0^1 (\xi(1)-\xi(q))\, d\ap (q).
	\end{equation*}
\end{lemma}
\textbf{Proof.}
Recall from \eqref{add:lem6:eq1}  and \eqref{SDEField} that $(\partial_x\Phi(s,X(s)))_{s\in[0,\qp]}$ reaches the optimal value of the stochastic optimal control problem for $\Phi(0,h)$, so 
\begin{align}
\Phi(0,h)&=\e \Phi(\qp ,X(\qp ))-\frac{1}{2}\int_0^{\qp }\xi''(s)\ap (s)\e\bigl[ \bigl(\partial_x\Phi(s,X(s))\bigr)^2\bigr]ds.
\label{Spart0}
\end{align}
Also, from \eqref{add:eq2},
\begin{align}
d\partial_x\Phi(s,X(s))&=\xi''(s)^{1/2}\partial_{xx}\Phi(s,X(s))dW(s),
\label{Spart1}
\\
d\partial_{xx}\Phi(s,X(s))&=-\xi''(s)\ap (s)\partial_{xx}\Phi(s,X(s))^2ds+\xi''(s)^{1/2}\partial_x^3\Phi(s,X(s))dW(s).
\label{Spart2}
\end{align}
Since $\ap(s)=1$ on $[\qp,1]$, the PDE can be explicitly solved up to $\qp$ via the Cole-Hopf transformation, 
\begin{align}
\label{add:eq-6}
\Phi(\qp , x) = \log\ch(x) + \frac{1}{2}\bigl(\xi'(1)-\xi'(\qp )\bigr).
\end{align}
As a result, $\myth(X(\qp ))=\partial_x\Phi(\qp ,X(\qp ))$ and
\begin{align*}
\e X(\qp )\myth(X(\qp ))
=&\,\, \e\Bigl[\Bigl(h+\int_0^{\qp }\xi''(s)\ap (s)\partial_x\Phi(s,X(s))ds\Bigr)\partial_x\Phi(\qp ,X(\qp ))\Bigr]
\\
&\,+\,\e\Bigl[\Bigl(\int_0^{\qp } \xi''(s)^{1/2}dW(s)\Bigr)\Bigl(\partial_x\Psi(0,X(0))+\int_0^{\qp }\xi''(s)^{1/2}\partial_{xx}\Phi(s,X(s))dW(s)\Bigr)\Bigr]
\\
=&\,\, h\e \myth(X(\qp ))+\int_0^{\qp }\xi''(s)\ap (s)\e\bigl[\bigl(\partial_x\Phi(s,X(s))\bigr)^2\bigr]ds
\\
&\,+\,\int_0^{\qp }\xi''(s)\e\bigl[\partial_{xx}\Phi(s,X(s))\bigr]ds,
\end{align*}
where the first and second equalities used \eqref{Spart1}. To handle the last two terms of the last equation, note that the equation (\ref{Spart1}) together with Ito's isometry gives
$$
\frac{d}{ds} \e\bigl[ \bigl(\partial_x\Phi(s,X(s)) \bigr)^2\bigr]
=
\xi''(s)\e\bigl[ \bigl(\partial_{xx}\Phi(s,X(s))\bigr)^2 \bigr].
$$
This and integration by parts give
\begin{align*}
&\int_0^{\qp }\xi''(s)\ap (s)\e\bigl[ \bigl(\partial_x\Phi(s,X(s)) \bigr)^2\bigr]ds\\
=&\,\, \xi'(s)\ap (s)\e\bigl[ \bigl(\partial_x\Phi(s,X(s)) \bigr)^2\bigr]\Big|_0^{\qp }
-\int_0^{\qp }\xi'(s)\e\bigl[ \bigl(\partial_x\Phi(s,X(s)) \bigr)^2\bigr] d\ap 
\\
&\,-\,\int_0^{\qp }\xi'(s)\xi''(s)\ap (s)\e\bigl[ \bigl(\partial_{xx}\Phi(s,X(s))\bigr)^2 \bigr]ds.
\end{align*}
On the other hand, using  (\ref{Spart2}) leads to
\begin{align*}
&\int_0^{\qp }\xi''(s)\e\bigl[\partial_{xx}\Phi(s,X(s))\bigr]ds\\
&=\xi'(s)\e\bigl[\partial_{xx}\Phi(s,X(s))\bigr]\Big|_{0}^{\qp }-\int_{0}^{\qp }\xi'(s)\frac{d}{ds}\e\bigl[\partial_{xx}\Phi(s,X(s))\bigr]ds\\
&=\xi'(\qp)\e\bigl[\partial_{xx}\Phi(\qp,X(\qp))\bigr]+\int_0^{\qp }\xi'(s)\xi''(s)\ap (s)\e\bigl[\bigl(\partial_{xx}\Phi(s,X(s))\bigr)^2\bigr]ds.
\end{align*}
Combining the above together yields
\begin{align*}
\e X(\qp )\myth(X(\qp ))&=h\e \myth(X(\qp ))+\xi'(\qp)\ap (\qp)\e\bigl[ \bigl(\partial_x\Phi(\qp,X(\qp)) \bigr)^2\bigr]\\
&-\int_0^{\qp }\xi'(s)\e\bigl[ \bigl(\partial_x\Phi(s,X(s)) \bigr)^2\bigr] d\ap+\xi'(\qp)\e\bigl[\partial_{xx}\Phi(\qp,X(\qp))\bigr].
\end{align*}
Now note that from the optimality of $\ap,$ $\e \partial_x\Phi(s,X(s))^2=s$ for all $s$ in the support of $\ap$ and that $\partial_x\Phi(\qp,x)=\myth(x)$ and $\partial_{xx}\Phi(\qp,x)=1-\myth^2(x)$. These and the above equation together with integration by part deduce that
\begin{align}
\begin{split}
\label{add:eq-8}
\e X(\qp )\myth(X(\qp ))&=h\e \myth(X(\qp ))+\xi'(\qp)\qp\ap (\qp)-\int_0^{\qp }\xi'(s)s d\ap+\xi'(\qp)(1-\qp)\\
&=h\e \myth(X(\qp ))+\int_0^{\qp }\xi''(s)\ap (s)sds+\int_0^{\qp }\xi'(s)\ap (s)ds+\xi'(\qp )(1-\qp ).
\end{split}
\end{align}
Next, recall that, by (\ref{Spart0}), \eqref{add:eq-6}, and the optimality of the Parisi measure,
\begin{align*}
\Phi(0,h)&=\e \log \ch( X(\qp ))+\frac{1}{2}(\xi'(1)-\xi'(\qp ))-\frac{1}{2}\int_0^{\qp }\xi''(s)\ap (s)sds.
\end{align*}
From this, since
\begin{align*}
\int_0^1\xi''(s)\ap (s)sds&=\int_0^{\qp }\xi''(s)\ap (s)sds+\int_{\qp }^1\xi''(s)sds\\
&=\int_0^{\qp }\xi''(s)\ap (s)sds+\xi'(s)s\Big|_{\qp }^1-\int_{\qp }^1\xi'(s)ds\\
&=\int_0^{\qp }\xi''(s)\ap (s)sds+\xi'(1)-\xi'(\qp )\qp -(\xi(1)-\xi(\qp )),
\end{align*}
the Parisi formula can be written as
\begin{align*}
\mathcal{P}(\ap )&:=\log 2+\Phi(0,h) - \frac{1}{2}\int_0^1\xi''(s)\ap (s)sds\\
&=\e \log 2\ch(X(\qp ))-\int_0^{\qp }\xi''(s)\ap (s)sds+C(\qp).
\end{align*}
Combining this with \eqref{add:eq-8} gives
\begin{align*}
\Delta&=\Bigl(\int_0^{\qp }\xi''(s)\ap (s)sds+\int_0^{\qp }\xi'(s)\ap (s)ds+\xi'(\qp )(1-\qp )\Bigr)\\
&\quad-\Bigl(\mathcal{P}(\ap )-C(\qp)+\int_0^{\qp }\xi''(s)\ap (s)sds\Bigr)\\
&=-\mathcal{P}(\ap )+C(\qp)+\int_0^{\qp }\xi'(s)\ap (s)ds+\xi'(\qp )(1-\qp )\\
&=-\mathcal{P}(\ap )-C(\qp)+\int_0^{\qp }\xi'(s)\ap (s)ds+\xi(1)-\xi(\qp).
\end{align*}
Finally, our proof is finished by substituting the following into the above equation,
\begin{align*}
\int_0^1 (\xi(1)-\xi(q))\, d\ap (q)&=\int_0^1\int_q^1\xi'(s)dsd\ap (q)\\
&=\int_0^1\xi'(s)\int_0^sd\ap (q)ds\\
&=\int_0^1 \xi'(s)\ap (s)ds\\
&=\int_0^{\qp }\xi'(s)\ap (s)ds+\xi(1)-\xi(\qp ).
\end{align*}
\qed

\subsection{Proof of Theorem \ref{ThHTtapLT}}

Recall $M_2(\sigma)$ from \eqref{SecTYk}.  Let $\qe$ be the EA order parameter associated to $m(\sigma).$ By definition, $\qe=M_2(\sigma).$ 
From \eqref{ad:eq-6} and \eqref{add:eq-7} with $k=2$ and the optimality of the Parisi measure,
\begin{align*}
\lim_{\varepsilon\downarrow 0}\lim_{N\rightarrow\infty}\e \bigl\la\bigl(C(\qe)-C(\qp)\bigr)^2\bigr\ra&\leq \frac{\xi''(1)}{2}\lim_{\varepsilon\downarrow 0}\lim_{N\rightarrow\infty}\e \bigl\la\bigl(\qe-\qp\bigr)^2\bigr\ra\\
&=\frac{\xi''(1)}{2}\lim_{\varepsilon\downarrow 0}\lim_{N\rightarrow\infty}\e \bigl\la\bigl(M_2(\sigma)-\qp\bigr)^2\bigr\ra=0.
\end{align*}
Write
\begin{align*}
F_N^{\tap}(m(\sigma))&=\frac{X_N(m(\sigma))}{N}-Y_N(\sigma)+C(M_2(\sigma)).
\end{align*}
From this expression, the above limit together with Lemmas \ref{LemFT}, \ref{LemST}, and \ref{lem6} implies that under the second moment with respect to the measure $\e\la \cdot\ra,$
 $F_N^{\tap}(m(\sigma))$ is equal to 
 \begin{align*}
\Bigl( \xi'(\qp)-\theta(\qp)-\int_0^1\xi(q)d\ap\Bigr)- \Delta+C(\qp)=\mathcal{P}(\ap)
 \end{align*}
 as $N$ tends to infinity and then $\varepsilon\downarrow 0.$ Finally, since $F_N$ converges to $\mathcal{P}(\ap)$ a.s. The dominated convergence theorem implies \eqref{ThHTtapLT:eq1}. This finishes our proof.

\bibliography{TAPFE}
\bibliographystyle{plain}
\end{document}